\documentclass[MSNbibl,number,citesort,dvips]{arxbj}
\usepackage{upgreek}
\usepackage{graphicx}

\aid{0}
\volume{20}
\issue{4}
\pubyear{2014}
\firstpage{2102}
\lastpage{2130}
\doi{10.3150/13-BEJ552} 

\makeatletter

\newcommand{\rright}{\right}
\newcommand{\lleft}{\left}
\newcommand{\rrVert}{\Vert}
\newcommand{\rrvert}{\vert}
\newcommand{\llVert}{\Vert}
\newcommand{\llvert}{\vert}

\newcommand{\idotsint}{\int\cdots\int}
\newcommand{\mathbbmathbf}[1]{\mathbf{#1}}
\newcommand{\trup}[2]{{#1}/{#2}}

\newtheorem{theorem}{Theorem}[section]
\newtheorem{corollary}[theorem]{Corollary}
\newtheorem{lemma}[theorem]{Lemma}
\newremark{Remark}{Remark}
\newremark{Lows}{Extreme significance levels}
\newremark{Deviationf}{Deviation from standard assumptions}
\newremark{Testpower}{Test power}
\newremark{Errorc}{Error-control quantities}
\newremark{Smalls}{Small sample sizes}
\newremark{Exacta}{Exact algebraic expression}
\newremark{Numericali}{Numerical integration using quadratures}
\newremark{Numericalint}{Numerical integration using Monte Carlo methods}
\newremark{Simulations}{Simulations}
\newremark{Estimation}{Estimation}
\newproclaim{Definition}[theorem]{Definition}
\makeatother

\begin{document}
\begin{frontmatter}

\title{Tail approximations for the Student $t$-, $F$-, and Welch
statistics for non-normal and not necessarily i.i.d. random variables}
\runtitle{Extremes of some popular statistical tests under non-normality}

\begin{aug}
\author{\inits{D.}\fnms{Dmitrii} \snm{Zholud}\corref{}\ead[label=e1]{dmitrii@chalmers.se}}
\address{Department of Mathematical Statistics,
Chalmers University of Technology and University of G\"{o}teborg,
SE-412 96 Gothenburg, Sweden. \printead{e1}}
\end{aug}

\received{\smonth{6} \syear{2012}}
\revised{\smonth{7} \syear{2013}}

%
\begin{abstract}
Let $T$ be the Student one- or two-sample $t$-, $F$-, or Welch
statistic.
Now release the underlying assumptions of normality, independence
and identical distribution
and consider a more general case where one only assumes that the
vector of data has a continuous joint density.
We determine asymptotic expressions for $\mathbf{P}(T>u)$ as $u\to
\infty$
for this case.
The approximations are particularly accurate for small sample sizes
and may be used, for example, in the analysis of
High-Throughput Screening experiments, where the number of replicates
can be as low as two to five
and often extreme  significance levels are used.
We give numerous examples and complement our results by an
investigation of the convergence speed -- both theoretically,
by deriving exact bounds for absolute and relative errors, and
by means of a simulation study.
\end{abstract}

%
\begin{keyword}
\kwd{dependent random variables}
\kwd{$F$-test}
\kwd{high-throughput screening}
\kwd{non-homogeneous data}
\kwd{non-normal population distribution}
\kwd{outliers}
\kwd{small sample size}
\kwd{Student's one- and two-sample $t$-statistics}
\kwd{systematic effects}
\kwd{test power}
\kwd{Welch statistic}
\end{keyword}

\end{frontmatter}

\section{Introduction}\vspace*{-2pt}

This article extends early results of Bradley
\cite{Bradley1952a} and Hotelling \cite{Hotelling1961} on
the tails of the distributions
of some popular and much used test statistics. We quantify the effect
of non-normality, dependence,
and non-homogeneity of data on the tails of the distribution of the
Student one- and two-sample $t$-, $F$- and Welch statistics. Our
approximations are valid for samples of any size,
but are most useful for very small sample sizes, for example, when
standard central limit theorem-based approximations are
inapplicable.\vspace*{-2pt}

\subsection{Problem statement and main result}

Let $\mathbf{X}\in\mathbb{R}^n$, $n\geq2$, be a random vector and
$T=T_n(\mathbf{X})$ be
(i) the Student one-sample $t$-test statistic; or (ii) the Student
two-sample $t$-test statistic; or
(iii) the $F$-test statistic for comparison of variances\vadjust{\goodbreak}
(in fact the $F$-test results apply also to one-way ANOVA, factorial designs,
a lack-of-fit sum of squares test, and an $F$-test for comparison of
two nested linear models).

In this paper, we study the asymptotic behavior of the tail
distribution of $T$ for small and fixed sample sizes.
Let $g_0(\mathbf{x})$ be the true joint density of $\mathbf{X}$ under
$H_0$ and $g_1(\mathbf{x})$ be the density under the alternative $H_1$.
Define $\mathcal{G}$ as a set of continuous densities that satisfy the
regularity constraints
of Theorems~\ref{OST_First_Order_Theorem}, \ref{TST_First_Order_Theorem},
or~\ref{FFirstOrderTheorem} for the
three test statistics accordingly.
Our main result is
the following theorem.\vspace*{-2pt}
%
%
\begin{theorem}
\label{IntroductionMainTheorem}
For any fixed value of $n$ and each of the three choices of $T$, there
exists a functional $K\dvtx\mathcal{G}\rightarrow\mathbb{R}^{+}$,
such that for all $g_0,g_1\in\mathcal{G}$ the limit expression
%
%
\begin{equation}
\label{Introduction_Main_result} \frac{\mathbf{P}({T>u|H_1})}{\mathbf{P}({T>u|H_0})}=\frac
{K_{g_1}}{K_{g_0}}+\mathrm{o}(1) \qquad\mbox{as }
u\to\infty
\end{equation}
holds with constants $0<K_{g_0}=K(g_0)<\infty$ and
$0<K_{g_1}=K(g_1)<\infty$.
The exact expressions for $K(g)$ are given in (\ref{OST_Def_Kg}),
(\ref{TST_Def_Kg}) and (\ref{F_Def_Kg}) for the three choices of the
test statistic~$T$.\vspace*{-2pt}
\end{theorem}

%
\begin{Remark}\label{rem1}
Standard assumption in the use of any of the test
statistics described above is that $\mathbf{X}\sim\operatorname
{MVN}(\mathbf
{0},\sigma^2\mathbf{1}_n)$, where $\operatorname{MVN}(\boldsymbol
\mu,\boldsymbol
\Sigma)$ denote the multivariate normal distribution with mean vector
$\boldsymbol\mu$ and covariance matrix $\boldsymbol\Sigma$.
It is easy to check that $\operatorname{MVN}(\mathbf{0},\sigma
^2\mathbf{1}_n)\in
\mathcal{G}$ and that $K(\operatorname{MVN}(\mathbf{0},\sigma
^2\mathbf{1}_n))=1$.\vspace*{-2pt}
\end{Remark}

Further remarks on Theorem~\ref{IntroductionMainTheorem} are given
in \hyperref[SectionSupplementaryMaterials]{Supplementary
Materials}, see \cite{supp}.\vspace*{-2pt}

\subsection{Motivation and applications}\label{sec1.2}\vspace*{-2pt}

The questions addressed in this article have gained significant new importance
through the explosive increase of High-Throughput Screening (HTS)
experiments, where the number of replicates
is often small, but instead thousands or millions of tests are
performed, at extremely high significance levels.
Studying extreme tails of test statistics under deviation from standard
assumptions is crucial in HTS because of the following factors:\vspace*{-2pt}

\begin{Lows*}
HTS uses many thousands or even
millions of biochemical, genetic or pharmacological tests.
In order to get a reasonable number of rejections, the significance
level of the tests is often very small,
say, 0.001 or lower, and it is the extreme tails of the distribution of
test statistics which are important.\vspace*{-2pt}
\end{Lows*}

\begin{Deviationf*}
HTS assays are often
subject to numerous systematic and spatial effects and to large number
of preprocessing steps.
The resulting data may become dependent, non-normal, or
non-homogeneous, yet common test statistics
such as one- and two-sample $t$-tests are still routinely computed
under standard assumptions.\vspace*{-2pt}
\end{Deviationf*}

\begin{Testpower*}
It is even less likely that the data follows any
standard distribution under the alternative hypothesis.
By quantifying the tail behavior of a test statistic under arbitrary
distributional assumptions, one can get more realistic estimators for
the test power.\looseness=-1\vadjust{\goodbreak}
\end{Testpower*}

\begin{Errorc*}
Given the scale of HTS experiments
and necessity to make even larger investments into further research on
positives detected through a HTS study, it is important to have
realistic picture of the accuracy of such experiments.
Consider, for example, estimation of $\mathit{pFDR}$, the positive False
Discovery Rate, see Storey \cite{Storey2002,Storey2003,Storey2004}. As of
now, estimates of $\mathit{pFDR}$ are obtained under the assumption
that the
true null distribution equals the theoretical one, and this may lead to
wrong decisions.
In most cases, however, a sample from the null distribution can be
obtained by conducting a separate experiment.
One can then model the tail distribution of the test statistic, and
apply for example, methods of Rootz\'{e}n and Zholud \cite{Rootzen2013},
which account for
deviations from the theoretical null distribution.
\end{Errorc*}

\begin{Smalls*} Due to economical constraints, numbers of
replicates in an individual experiment in HTS are as small as two to five,
which makes large sample normal approximations inapplicable. Even for
moderate sample sizes, CLT-based approximations are not accurate in the
tails and
better approximations, such as those presented in this paper, are needed.
\end{Smalls*}

\begin{figure}

\includegraphics{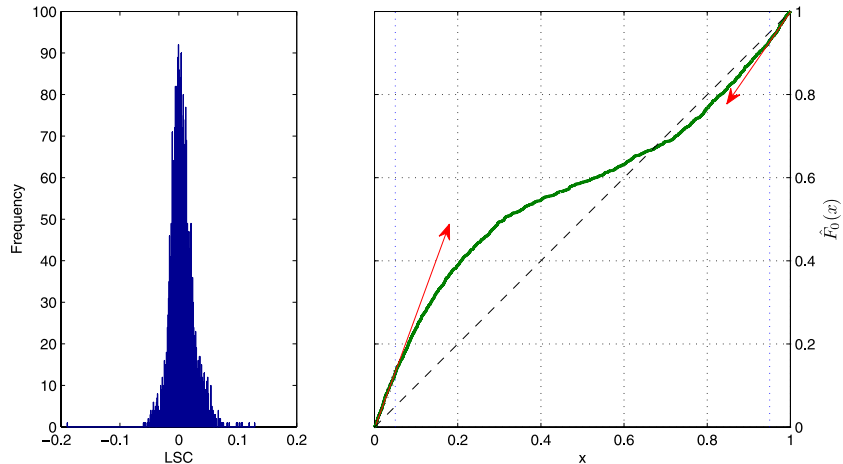}

\caption{The wild type data set. \textit{Left}: Histogram of $3456$
LSC values from the wildtype dataset.
\textit{Right}: Empirical CDF of 1728 $p$-values obtained from
one-sample $t$-test for pairs of LSC values.}\label{fig:BioscreenECDF}
\end{figure}

We now consider a HTS experiment which was the motivation for the
present paper.
Left panel of Figure~\ref{fig:BioscreenECDF} shows measured values of
the \textit{Logarithmic Strain Coefficient} (LSC) of
the wildtype cells in a Bioscreen array experiment in yeast genome
screening studies, see Warringer and Blomberg \cite{Warringer2003a}
and Warringer \textit{et al.} \cite{Warringer2003b}.

%
%

The null hypothesis was that LSC of a wildtype yeast cell had normal
distribution with mean zero and unknown variance.
The experiment was made for quality control purposes, hence no
treatment has been applied and the null hypothesis of mean zero was
known to be true.

The histogram of the LSC values was skewed, see Figure~\ref
{fig:BioscreenECDF},
and we therefore plotted the empirical cumulative distribution function
(CDF) of the $1728$ $p$-values obtained from the LSC values.
As expected, see the right panel in Figure~\ref{fig:BioscreenECDF},
the distribution of the $p$-values was different from the theoretical
uniform distribution.

Note, however, that both lower and upper tails of the plot approach
straight lines, as indicated by the two arrows. This was in fact the
starting point of the present article, and it later
followed that such tail behavior is justified by
Theorem~\ref{IntroductionMainTheorem}, see
\hyperref[SectionSupplementaryMaterials]{Supplementary Materials}.

In practical applications, one needs to be able to compute or estimate
the constant $K_{g}$. This can be done in a variety of ways.

\begin{Exacta*}
For the case when components of
$\mathbf{X}$ are i.i.d. random variables, constant $K_g$ can be
obtained directly from (\ref{OST_Def_Kg}), (\ref{TST_Def_Kg}) and
(\ref{F_Def_Kg}) for the three choices of the test statistic
$T$\vadjust{\goodbreak}
accordingly. We give numerous examples through Sections~\ref
{Section_One_Sample}--\ref{Section_F},
and \hyperref[SectionSupplementaryMaterials]{Supplementary Materials}
provides Wolfram {M}athematica \cite{Mathematica2010} code to
compute $K_g$ for even more complicated cases,
like, for example, Multivariate Normal case with $g\sim\operatorname
{MVN}(\mathbf
{0},\boldsymbol{\Sigma})$.
\end{Exacta*}

\begin{Numericali*}
For an arbitrary
multivariate density $g(\mathbf{x})$ and Student one- and two-sample
$t$-statistics, or $F$-statistics with low degrees of freedom, $K_g$
can be computed from (\ref{OST_Def_Kg}), (\ref{TST_Def_Kg}) or
(\ref{F_Def_Kg}) using adaptive Simpson or Lobatto quadratures. We provide
the corresponding {MATLAB} \cite{Matlab2010} scripts in
\hyperref[SectionSupplementaryMaterials]{Supplementary Materials}.
\end{Numericali*}

\begin{Numericalint*}
For an
$F$-statistic with the denominator that has more than two degrees of
freedom, $K_g$ can be computed numerically using Monte Carlo
integration, see
\hyperref[SectionSupplementaryMaterials]{Supplementary Materials}.
Monte Carlo
methods are applicable to the case described above as well.
\end{Numericalint*}

\begin{Simulations*}
The distribution tail of $T$ can be estimated
using simulations, see, for example, Section~\ref{SectionSimulation}.
In the current paper we used ``brute-force'' approach, but importance
sampling techniques can be applied quite generally as well.
\end{Simulations*}

\begin{Estimation*}
If $g(\mathbf{x})$ is unknown but one instead
has a sample from $g$, then $K_g$ can be estimated as a slope of the
graph of the CDF of the corresponding $p$-values in the origin of zero.
In the yeast genome screening experiment, for example, $K_g$
approximately equals the slope of the red arrow -- theoretical
justification of this fact is given in
\hyperref[SectionSupplementaryMaterials]{Supplementary Materials},
and the
estimation technique is similar to the Peak-Over-Threshold (POT) method
in Extreme Value Theory, see, for example, the SmartTail
software at
\href{http://www.smarttail.se}{www.smarttail.se} \cite{SmartTail} and
further examples in Rootz\'{e}n and Zholud \cite{Rootzen2013}.
\end{Estimation*}

Finally, the existence of $K_g$ and its importance for questioning the
logic behind some multiple testing procedures is discussed in
Zholud \cite{Zholud2011a}, Part~I, Section~3.

\subsection{Literature review}

There is enormous amount of literature on the behavior of the Student
one- and two-sample $t$- and $F$-statistics under deviations from the
standard assumptions.
The overwhelming part of this literature is focused on normal
approximations, that is, when $n\to\infty$. These are large sample
approximations though,
and are irrelevant to the topic of the present article.

For small and moderate sample sizes, one would typically use Edgeworth
expansion, see, for example, Field and Ronchetti \cite{Field1990},
Hall \cite{Hall1987} and
Gaen \cite{Gaen1949,Gaen1950}, or
{saddlepoint approximations}, see, for
example, Zhou and Jing \cite{Zhou2006},
Jing \textit{et al.} \cite{Jing2004} and Daniels and Young \cite
{Daniels1991}.
Edgeworth expansion improves the normal approximation but is still
inaccurate in the tails.
Saddlepoint approximations, on the other hand, can be very accurate in
the tails, see, for example, Jing \textit{et al.} \cite{Jing2004},
but the latter
statement is based on purely empirical evidence and the asymptotic
behavior of these approximations as $u\to\infty$ is not well studied.
Furthermore, in practice one would require exact parametric form of the
population density, and the use of saddle point approximations in
statistical inference is questionable.\looseness=-1

As for the approximations considered in this article, that is, when $n$
is small and $u\to\infty$, the existing literature is very limited.
This presumably can be explained by the fact that
situations where one would need to test at significance levels of
$10^{-3}$ and lower never arose, until present times. We focus on
the most relevant works by Bradley \cite{Bradley1952a,Bradley1952b} and
Hotelling \cite{Hotelling1961}.

Bradley covers the Student one-sample $t$-statistic
for i.i.d. non-normal observations, and also makes a somewhat less
complete study of the corresponding cases for the Student two-sample
$t$-test and the $F$-test of equality of variances.
Bradley \cite{Bradley1952b} derives the constant
$K_g$ from geometrical
considerations, but does not state any assumptions on the underlying
population density which ensure that the approximations hold.
Bradley \cite{Bradley1952a}, on the other hand,
gives assumptions on the population
density, but
these assumptions are insufficient, see Section~\ref{Section_Understanding_The_Technical_Restriction}.

Hotelling \cite{Hotelling1961} studies the Student
one-sample $t$-test for an
``arbitrary'' joint density of $\mathbf{X}$.
Hotelling derives the constant $K_g$ assuming that
the limit in the left-hand side of
(\ref{Introduction_Main_result}) exists and that the function
\[
D_n(\xi)=\int_0^{\infty}r^{n-1}g(r
\xi_1,\ldots,r\xi_n)\,\mathrm{d}r
\]
is continuous for both densities $g_0$ and $g_1$.
When it comes to the examples, however, the existence of the limit in
(\ref{Introduction_Main_result}) is taken for granted and the
assumption of continuity of $D_n(\xi)$ is never verified.

Finally, a more detailed literature review that covers other approaches
and meritable scientific works is given in
\hyperref[SectionSupplementaryMaterials]{Supplementary Materials}.

The structure of this paper is as follows:
Sections~\ref{Section_One_Sample}--\ref{Section_F} contain main
theorems and examples;
Section~\ref{Section_Second_Order_Approximations} addresses the
convergence speed and higher order expansions;
Section~\ref{SectionSimulation} presents a simulation study.
\hyperref[Section_FurtherRemarks]{Appendix A} includes the key lemma
used in the proofs, in Section~\ref{Section_Shrinking_Balls}, and a discussion on the regularity
conditions, in Section~\ref{Section_Understanding_The_Technical_Restriction};
\hyperref[SectionFiguresandTables]{Appendix B} contains figures
from the
simulation study; and, finally, follows a brief summary of the
\hyperref[SectionSupplementaryMaterials]{Supplementary Materials}
that are available online.
%
\section{One-sample $t$-statistic}%
\label{Section_One_Sample}%

Let
$\mathbf{X}=(X_1,X_2,\ldots,X_n)$, $n\geq2$,
be a random vector that has a joint density $g$
and define
\begingroup
\abovedisplayskip=7.5pt
\belowdisplayskip=7.5pt
\[
T = \sqrt{n} ( \overline{\mathbf{X}} / {S}),
\]
where $\overline{\mathbf{X}}$ and $S^2$ are the sample mean and the
sample variance
of the vector $\mathbf{X}$.
Introduce the unit vector
$
\mathbbmathbf{I}= (1/{\sqrt{n}},1/{\sqrt{n}},\ldots
,1/{\sqrt{n}} )$,
and assume that
%
%
\begin{equation}
\label{OST_Condition_On_Density_Positive_on_Diagonal} g(x\mathbbmathbf{I})>0 \qquad\mbox{for some } x\geq0
\end{equation}
and that
%
%
\begin{equation}
\label{OST_Condition_On_Density} \int_{0}^{\infty} r^{n-1}
\mathop{\sup_{\|{\boldsymbol{\xi}}\|<\varepsilon,}}_{{\boldsymbol
{\xi}}\in
L^{\bot}} g \bigl( r (
\mathbbmathbf{I} + {\boldsymbol{\xi}} ) \bigr) \,\mathrm{d}r < \infty
\end{equation}
for some $\varepsilon>0$, where $L$ is the linear subspace of $\mathbb
{R}^{n}$ spanned by the vector $\mathbbmathbf{I}$ and $L^\bot$ is
its orthogonal complement.
Finally, introduce the constant
%
%
\begin{equation}
\label{OST_Def_Kg} K_g= 2 \frac{
\uppi^{\trup{n}{2}}
} {
\Gamma(\trup{n}{2})
} \int_{0}^{\infty}
r^{n-1} g (r\mathbbmathbf{I} ) \,\mathrm{d}r.
\end{equation}

%
\begin{theorem}
\label{OST_First_Order_Theorem}
If $g$ is continuous and satisfies
(\ref{OST_Condition_On_Density_Positive_on_Diagonal}) and
(\ref{OST_Condition_On_Density}), then
%
%
\begin{equation}
\label{OST_First_Order} \frac{\mathbf{P}({T>u})}{t_{n-1}(u)}=K_g+\mathrm{o}(1) \qquad\mbox {as
} u\to\infty,
\end{equation}
where $t_{n-1}(u)$ is the tail of the $t$-distribution with $n-1$
degrees of freedom and $0<K_g=\allowbreak  K(g)<\infty$.
\end{theorem}
\begin{pf}
We use several variable changes to transform the right-hand side of
\[
\mathbf{P}({T>u})=\int_{D_1} g(\mathbf{x}) \,\mathrm{d}
\mathbf{x},
\]
where
$D_1= \{\mathbf{x}\dvtx T>u  \}$
and
$\mathrm{d}\mathbf{x}$
is the notation for
$\mathrm{d}x_1\,\mathrm{d}x_2\cdots\,\mathrm{d}x_{n}$, to the form
treated in Corollary~\ref{IOFOSB_Base_Lemma_Corollary}.
Let
$\mathbf{e}_{\mathbf{1}},\mathbf{e}_{\mathbf{2}},\ldots,
\mathbf{e}_{\mathbf{n}}$
be the standard basis in
$\mathbb{R}^{n}$
and $A$ be an orthogonal linear operator which satisfies
%
%
\begin{equation}
\label{OST_Def_A} A\mathbf{e}_{\mathbf{n}}=\mathbbmathbf{I}.
\end{equation}
Setting $\mathbf{x}=A\mathbf{y}$ we have that $\overline{\mathbf
{X}}=y_{n}/\sqrt{n}$ and $S^2=\sum_{i=1}^{n-1}y^2_i/(n-1)$, and hence
\[
\mathbf{P}({T>u})=\int_{D_2} g(A\mathbf{y}) \,\mathrm{d}
\mathbf{y},
\]
where
\[
D_2=
\biggl\{
\mathbf{y}\dvtx
\frac{
y_{n}
} {
\sqrt{\trup{1}{(n-1)}\sum_{i=1}^{n-1}y^2_i}
}
>u
\biggr\}.\vadjust{\goodbreak}
\]

Next, introducing new variables $y_i=(n-1)^{1/2} r t_i$ for $i\leq n-1$ and
$y_{n}=r$, $r>0$, applying Fubini's theorem, and recalling
(\ref{OST_Def_A}) we get
%
%
\begin{equation}
\label{OST_Iterated_int} \mathbf{P}({T>u}) = \mathop{\idotsint}_{\sum t_i^2<u^{-2}} G(
\mathbf{t}) \,\mathrm{d}\mathbf{t},
\end{equation}
where
\[
G(\mathbf{t})= (n-1)^{\trup{(n-1)}{2}} \int_0^\infty
r^{n-1} g \bigl( r \bigl( \mathbbmathbf{I} + A\mathbf{v}(\mathbf{t})
\bigr) \bigr) \,\mathrm{d}r,
\]
and
\[
\mathbf{v}(\mathbf{t}) = (n-1)^{1/2} ( t_1,t_2,
\ldots,t_{n-1},0 ).
\]
Continuity of $g$ and (\ref{OST_Condition_On_Density}) ensure that $G$
is continuous at zero, by
the dominated convergence theorem, and Corollary~\ref{IOFOSB_Base_Lemma_Corollary}
completes the proof.
\end{pf}\endgroup

Assumption (\ref{OST_Condition_On_Density_Positive_on_Diagonal})
ensures that $K_g>0$ and the condition
(\ref{OST_Condition_On_Density}) holds if, for example,
$K_g<\infty$ and $g$ is continuous and has the asymptotic monotonicity
property, see
Lemma~\ref{Asmptotic_Monotonicity_Criterium}.

Now consider the case when one of the assumptions
(\ref{OST_Condition_On_Density}) or
(\ref{OST_Condition_On_Density_Positive_on_Diagonal})
is violated.
If (\ref{OST_Condition_On_Density}) holds and
(\ref{OST_Condition_On_Density_Positive_on_Diagonal})
is violated,
then (\ref{OST_First_Order}) holds with $K_g=0$, that is, the right
tail of the distribution of $T$ is ``strictly lighter'' than
$t_{n-1}(u)$, the tail of the $t$-distribution with $n-1$ degrees of freedom.
If, instead, (\ref{OST_Condition_On_Density_Positive_on_Diagonal})
holds and (\ref{OST_Condition_On_Density}) is violated, then,
Theorem~\ref{Relaxing_Assumptions_LimInf_Theorem}
shows that
the right tail of the distribution of $T$ is ``at least as heavy'' as
$t_{n-1}(u)$, provided $K_g<\infty$, and ``strictly heavier''
than $t_{n-1}(u)$ if $K_g=\infty$.

We next consider two important corollaries -- one concerning dependent
Gaussian vectors, and another one that addresses the non-normal i.i.d. case.
%
%
\begin{corollary}[(Gaussian zero-mean case)]\label{OST_GaussianCorollary}
If $\mathbf{X}\sim\operatorname{MVN}(\boldsymbol{0},\boldsymbol
{\Sigma})$, where
$\boldsymbol{\Sigma}$ is a strictly positive-definite covariance
matrix, then (\ref{OST_First_Order}) holds with
\[
K_g=\frac{
(\mathbbmathbf{I} {\boldsymbol{\Sigma}} \mathbbmathbf
{I}^T )^{n/2}
} {
|\boldsymbol{\Sigma}|^{1/2}
}.
\]
\end{corollary}
\begin{pf}
Deriving the expression for $K_g$ in (\ref{OST_Def_Kg}) is
straightforward. Note that $K_g<\infty$
since ${\boldsymbol{\Sigma}}$ is non-degenerate and
$\operatorname{MVN}(\boldsymbol
{0},\boldsymbol{\Sigma})$ has
the asymptotic monotonicity property defined in
Section~\ref{Section_Understanding_The_Technical_Restriction}.
It then follows from Lemma~\ref{Asmptotic_Monotonicity_Criterium} that
the regularity constraint (\ref{OST_Condition_On_Density}) holds, and
so does~(\ref{OST_First_Order}).
\end{pf}

One possible application of Corollary~\ref{OST_GaussianCorollary} is
to correct for the effect of dependency when using test statistic $T$.
This is done by dividing the corresponding $p$-value by $K_g$.

Now consider the effect of non-normality. Assume that the elements
$X_i$ of the vector $\mathbf{X}$ are independent and
identically distributed and let $h(x)$ be their common marginal
density, so that $g(\mathbf{x})=h(x_1)h(x_2)\cdots h(x_n)$.
%
%
\begin{corollary}[(i.i.d. case)]\label{OST_Corollary_iid}
If $h(x)$ is continuous, and monotone on $[L,  \infty)$ for some
finite constant $L$, then (\ref{OST_First_Order}) holds with
\[
K_g= 2 \frac{
(\uppi n)^{\trup{n}{2}}
} {
\Gamma(\trup{n}{2})
} \int_{0}^{\infty}
r^{n-1} h (r )^n \,\mathrm{d}r<\infty.
\]
\end{corollary}
\begin{pf}
The monotonicity of $h(x)$ on $[L,\infty)$ implies that $g(\mathbf
{x})$ has the asymptotic monotonicity property,
see
Section~\ref{Section_Understanding_The_Technical_Restriction},
and the regularity assumption
(\ref{OST_Condition_On_Density})
hence follows from finiteness of $K_g$ and
Lemma~\ref{Asmptotic_Monotonicity_Criterium}.
The finiteness of $K_g$, in turn, follows if we show that $rh(r)\to0$
as $r\to\infty$.

Indeed, assume to the contrary that $\limsup rh(r)>0$. Then there
exists $\delta>0$ and a sequence $\{r_k\}_{k=0}^{\infty}$ with $r_0=L+1$
and such that $r_{k+1}>2r_k$ and $r_kh(r_k)>\delta$ for any $k>0$.
Now the monotonicity of $h(x)$ on $[L,\infty)$ gives
\[
\int_{L+1}^{\infty} h(r) \,\mathrm{d}r \geq \sum
_{k=1}^{\infty} (r_{k}-r_{k-1}
)h(r_k) > \delta\sum_{k=1}^{\infty}
\frac{r_{k}-r_{k-1}}{r_k} = \infty,
\]
contradicting that $h(x)$ is a density.
\end{pf}

The constants $K_g$ for some common densities $h(x)$ are given in
Table~\ref{tab:KgOSTIID}.

%
%
\begin{table}
\tablewidth=\textwidth
\tabcolsep=0pt
\caption{The constants $K_g$ for the i.i.d. case of the Student
one-sample $t$-test. Here
$\Gamma(x)$, $B(x)$ and $M(a,b,x)$ are the Gamma, Beta and Kummer
confluent hypergeometric function, see, for example,
Hayek~\cite{Hayek2001}}\label{tab:KgOSTIID}
\begin{tabular*}{\textwidth}{@{\extracolsep{\fill}}l@{}}
\hline
Normal with mean $\mu\ne0$ and standard deviation $\sigma>0$
\\
$M
(\frac{1-n}{2},
\frac{1}{2},
-\frac{n\mu^2}{2\sigma^2}
)
+
\frac{\mu}{\sigma}\frac{\sqrt{2n}\Gamma
(\trup{(1+n)}{2} )} {
\Gamma (\trup{n}{2} )}
M
(
1-\frac{n}{2},
\frac{3}{2},
-\frac{n\mu^2}{2\sigma^2}
)
$
\\[9pt]
Half-normal, and log-normal derived from a $N(\mu,\sigma^2)$
\\
$2^n$
and
$\frac{n^{\trup{(n-1)}{2}}\sqrt{\uppi}} {
2^{\trup{(n-3)}{2}}
\sigma^{n-1}\Gamma (\trup{n}{2} )}$
\\[9pt]
$\chi$ with $\nu>0$, and $\chi^2$ (and its inverse) with $\nu\geq
2$ d.f.
\\
$\frac{2^n\uppi^{n/2}\Gamma
(\trup{n \nu}{2} )}{n^{\trup{n}{2}(\nu-1)}
\Gamma (\trup{\nu}{2} )^{n}
\Gamma (\trup{n}{2} )}$
and
$\frac{2\uppi^{n/2}\Gamma (\trup{n \nu}{2} )} {
n^{(\trup{n}{2})(\nu-1)}\Gamma
(\trup{\nu}{2} )^{n}
\Gamma (\trup{n}{2} )}$
\\[9pt]
$F$ with $\mu>0$ and $\nu>0$ degrees of freedom
\\
$\frac{2(\uppi n)^{n/2}\Gamma
(\trup{\mu n}{2} )
\Gamma (\trup{\nu n}{2} )
\Gamma (\trup{(\mu+\nu)}{2} )^n} {
\Gamma (\trup{n}{2} )
[\Gamma (\trup{\mu}{2} )
\Gamma (\trup{\nu}{2} )
]^{n}
\Gamma (\trup{(\mu+\nu)}{2} n  )}$
\\[9pt]
$T$ with $\nu>0$ d.f. and Cauchy
\\
$\frac{n^{n/2}\Gamma (\trup{\nu n }{2} )}{\Gamma
(\trup{(\nu+1)n}{2} )}
(\frac{\Gamma (\trup{(\nu+1)}{2} )} {
\Gamma (\trup{\nu}{2} )} )^n$
and
$\frac{n^{n/2}}{2^{n-1}\uppi^{\trup{(n-1)}{2}}\Gamma (\trup
{(n+1)}{2} )}$
\\[9pt]
Beta with shape parameters $\alpha>1$ and $\beta>1$
\\
$\frac{2(\uppi n)^{n/2}\Gamma(\alpha n)
\Gamma (1+(\beta-1)n  )} {
B(\alpha,\beta)^{n}\Gamma (\trup{n}{2} )
\Gamma (1+(\alpha+\beta-1)n
)}$
\\[9pt]
Gamma (and its inverse) with shape $\alpha>1$
\\
$\frac{2n^{\trup{n}{2}(1-2\alpha)}
\uppi^{n/2}\Gamma(\alpha n)}{\Gamma(\alpha)^{n}
\Gamma (\trup{n}{2} )}$
\\[9pt]
Uniform on interval $[a,b]$, $b>0$
\\
$\frac{(\uppi n)^{\trup{n}{2}}} {
\Gamma (\trup{n}{2}+1 )}
\cases{
(\frac{b}{b-a} )^n &\quad$0\in[a, b]$,
\cr
\frac{b^n-a^n}{(b-a)^{n}} &\quad$[a,b]\subset[0,\infty)$}$
\\[14pt]
Centered exponential and exponential
\\
$\frac{2 (\trup{\uppi}{n} )^{n/2}\Gamma(n)} {
\mathrm{e}^{n}\Gamma (\trup{n}{2} )}$
and
$\frac{2 (\trup{\uppi}{n} )^{n/2}\Gamma(n)} {
\Gamma (\trup{n}{2} )}$
\\[9pt]
Maxwell, and Pareto with $k>0$ and scale $\alpha>0$
\\
$\frac{ (\trup{4}{n} )^n
\Gamma (\trup{3 n}{2} )} {
\Gamma (\trup{n}{2} )}$
and
$\frac{(\uppi n)^{n/2}\alpha^{n-1}} {
\Gamma (\trup{n}{2}+1 )}$
\\
\hline
\end{tabular*}
\end{table}

\section{Two-sample $t$-statistic}\label{Section_Two_Sample}%

In this section, we cover the Student two-sample $t$-statistic.
However, we first consider a more general case.
For $n_1\geq2$, $n_2\geq2$, set $n=n_1+n_2$ and let $\mathbf
{X}=(X_1,X_2,\ldots,X_{n})$
be a random vector that has a multivariate joint density $g$.
Further, let $S_1$ and $S_2$ be the sample variances
of the vectors $(X_1,X_2,\ldots,X_{n_1})$ and
$(X_{n_1+1},X_{n_1+2},\ldots,X_{n})$ and define
\[
T = \frac{
(\trup{1}{n_1}) \sum_{i=1}^{n_1}X_i - (\trup{1}{n_2})\sum_{i=n_1+1}^{n}X_i
} {
\sqrt{\alpha S_1^2 + \beta S_2^2}
},
\]
where $\alpha$ and $\beta$ are some positive constants (to be set later).
Next, define the two unit vectors
$
\mathbbmathbf{I}_{1}\!=\! (1/{\sqrt{n_1}},1/{\sqrt{n_1}},\ldots
,1/{\sqrt{n_1}},0,0,\ldots,0 )
$
and
$
\mathbbmathbf{I}_{2}\!=\! (0,0,\ldots,0,1/{\sqrt{n_2}},1/{\sqrt {n_2}},\ldots,\break 1/{\sqrt{n_2}} )$,
and let $\omega_0=\arccos (\sqrt{{n_2}/{n}} )$. We assume that
%
%
\begin{equation}
\label{TST_Condition_On_Density_Positive_on_the_plane} g \bigl( r \bigl( \cos(\omega-\omega_0)
\mathbbmathbf{I}_{1} + \sin(\omega-\omega_0)
\mathbbmathbf{I}_{2} \bigr) \bigr)>0
\end{equation}
for some $r\geq0$ and $\omega\in[-\uppi/2,  \uppi/2]$, and that for
some $\varepsilon>0$
%
%
\begin{eqnarray}
\label{TST_Condition_On_Density} &&\int_{-\uppi/2}^{\uppi/2} \cos(
\omega)^{n-2}\int_0^\infty
r^{n -1}
\nonumber
\\[-8pt]
\\[-8pt]
&&\phantom{\int_{-\uppi/2}^{\uppi/2} \cos(
\omega)^{n-2}\int_0^\infty}{}\times\mathop{\sup_{\|{\boldsymbol\xi}\|<\varepsilon
}}_{{\boldsymbol\xi
}\in L^\bot} g \bigl( r
\bigl( \cos(\omega-\omega_0)\mathbbmathbf{I}_{1} + \sin(
\omega-\omega_0)\mathbbmathbf{I}_{2} + {\boldsymbol{\xi}}
\bigr) \bigr) \,\mathrm{d}r \,\mathrm{d}\omega<\infty,
\nonumber
\end{eqnarray}
where $L$ is a linear subspace of $\mathbb{R}^{n}$ spanned by the
vectors $\mathbbmathbf{I}_{1}$ and $\mathbbmathbf{I}_{2}$,
and $L^\perp$ is its orthogonal complement. Next, define the constant
%
%
\begin{eqnarray}
\label{TST_Def_Kg} K_g&=&C(n_1,n_2,\alpha,
\beta) \int_{-\uppi/2}^{\uppi/2} \cos(\omega)^{n-2}
\nonumber
\\[-8pt]
\\[-8pt]
&&\phantom{C(n_1,n_2,\alpha,
\beta) \int_{-\uppi/2}^{\uppi/2}}{}\times \int_0^\infty r^{n -1} g
\bigl( r \bigl( \cos(\omega-\omega_0)\mathbbmathbf{I}_{1}
+ \sin(\omega-\omega_0)\mathbbmathbf{I}_{2} \bigr) \bigr)
\,\mathrm{d}r \,\mathrm{d}\omega,
\nonumber
\end{eqnarray}
where the constant $C(n_1,n_2,\alpha,\beta)$ is given by
\[
C(n_1,n_2,\alpha,\beta)= \frac{2
\uppi^{\trup{(n-1)}{2}}
(\trup{(n_1-1)}{\alpha} )^{\trup{(n_1-1)}{2}}
(\trup{(n_2-1)}{\beta} )^{\trup{(n_2-1)}{2}}
(\trup{1}{n_1}+\trup{1}{n_2} )^{\trup{(n-2)}{2}}
} {
\Gamma (\trup{(n-1)}{2} )
(n-2)^{\trup{(n-2)}{2}}
}.
\]

%
\begin{theorem}
\label{TST_First_Order_Theorem}
If $g$ is continuous and satisfies
(\ref{TST_Condition_On_Density_Positive_on_the_plane}) and
(\ref{TST_Condition_On_Density}), then
%
%
\begin{equation}
\label{TST_First_Order} \frac{\mathbf{P}({T>u})}{t_{n-2}(u)}=K_g+\mathrm{o}(1) \qquad\mbox {as
} u\to\infty,
\end{equation}
where $t_{n-2}(u)$ is the tail of the $t$-distribution with $n-2$
degrees of freedom and $0<K_g=\allowbreak  K(g)<\infty$.
\end{theorem}

\begin{pf}
The proof is similar to the proof of Theorem~\ref{OST_First_Order_Theorem}.
Let $A$ be an orthogonal linear operator such that
%
%
\begin{equation}
\label{TST_Def_A} A\mathbf{e}_{\mathbf{n}_{\mathbf{1}}}=\mathbbmathbf{I}_{1} \quad
\mbox{and}\quad A\mathbf{e}_{\mathbf{n}}=\mathbbmathbf{I}_{2}.
\end{equation}
Changing coordinate system $\mathbf{x}=A\mathbf{y}$ gives
\begin{eqnarray*}
\frac{1}{n_1}\sum_{i=1}^{n_1}X_i&=&y_{n_1}/
\sqrt{n_1}, \qquad \frac{1}{n_2}\sum
_{i=n_1+1}^{n}X_i=y_{n}/
\sqrt{n_2},
\\
S^2_1&=&\sum_{i=1}^{n_1-1}y^2_i/(n_1-1)
\quad \mbox{and} \quad S^2_2=\sum
_{i=n_1+1}^{n-1}y^2_i/(n_2-1)
\end{eqnarray*}
and therefore
\[
\mathbf{P}({T>u})=\int_{ \{\mathbf{x}: T>u  \}} g(\mathbf{x}) \,\mathrm{d}
\mathbf{x}=\int_{D} g(A\mathbf{y}) \,\mathrm{d}\mathbf{y},
\]
where
\[
D= \Biggl\{ \mathbf{y}\dvtx \biggl(\frac{1}{\sqrt{n_1}}y_{n_1}-
\frac{1}{\sqrt{n_2}}y_{n} \biggr) \bigg/ \Biggl( \frac{\alpha}{n_1-1}\sum
_{i=1}^{n_1-1}y^2_i+
\frac{\beta}{n_2-1}\sum_{i=n_1+1}^{n-1}y^2_i
\Biggr)^{1/2} >u  \Biggr\}.
\]
Next, define $c_1(\omega)$ and $c_2(\omega)$ by
\[
\frac{c_1(\omega)}{\sqrt{{1}/{n_1} + {1}/{n_2}}}= \sqrt{\frac{n_1-1}{\alpha}} \cos(\omega) \quad\mbox{and} \quad
\frac{c_2(\omega)}{\sqrt{{1}/{n_1} + {1}/{n_2}}}= \sqrt{\frac{n_2-1}{\beta}} \cos(\omega),
\]
and introduce new variables $t_1,t_2,\ldots,t_{n-2},r,\omega$ such that
\begin{eqnarray*}
y_i &=&r c_1(\omega) t_i \qquad\mbox{for
} i=1,2,\ldots,n_1-1,
\\
y_i &=&r c_2(\omega) t_{i-1} \qquad\mbox{for
} i=n_1+1,n_1+2,\ldots,n-1,
\\
y_{n_1}&=&r \cos(\omega-\omega_0) \quad\mbox{and}\quad
y_{n}=r \sin(\omega -\omega_0),\qquad r>0.
\end{eqnarray*}
The identity
$
\cos(\omega-\omega_0)/{\sqrt{n_1}}-\sin(\omega-\omega_0)/{\sqrt{n_2}}
=
\sqrt{{1}/{n_1}+{1}/{n_2}}\cos(\omega)$,
Fubini's theorem, and (\ref{TST_Def_A}) give
%
%
\begin{equation}
\label{TST_Iterated_int} \mathbf{P}({T>u}) = \mathop{\idotsint}_{\sum_{i=1}^{n-2} t_i^2<u^{-2}} G(
\mathbf{t}) \,\mathrm{d}\mathbf{t},
\end{equation}
where
\begin{eqnarray*}
G(\mathbf{t})&=& M \int_{-\uppi/2}^{\uppi/2} \cos(
\omega)^{n-2} \int_0^\infty
r^{n -1}
\\
&&\phantom{M \int_{-\uppi/2}^{\uppi/2} \cos(
\omega)^{n-2} \int_0^a}{}\times g \bigl( r \bigl( \cos(\omega-\omega_0)
\mathbbmathbf{I}_{1} + \sin(\omega-\omega_0)
\mathbbmathbf{I}_{2} + A\mathbf{v}(\mathbf{t},\omega-
\omega_0) \bigr) \bigr) \,\mathrm{d}r \,\mathrm{d}\omega
\end{eqnarray*}
with
\[
\mathbf{v}(\mathbf{t},\omega) = \bigl( c_1(\omega)t_1,
\ldots,c_1(\omega)t_{n_1-1} , 0 , c_2(
\omega)t_{n_1},\ldots,c_2(\omega)t_{n-2} , 0 \bigr)
\]
and
\[
M = \biggl(\frac{n_1-1}{\alpha} \biggr)^{\trup{(n_1-1)}{2}} \biggl(\frac{n_2-1}{\beta}
\biggr)^{\trup{(n_2-1)}{2}} \biggl(\frac{1}{n_1}+\frac{1}{n_2}
\biggr)^{\trup{(n-2)}{2}}.
\]
The finiteness of the integral in (\ref{TST_Condition_On_Density}) and
continuity of $g$ imply the
continuity of $G$ at zero by the dominated convergence theorem, and
Corollary~\ref{IOFOSB_Base_Lemma_Corollary}
gives the asymptotic
expression (\ref{TST_First_Order}) with the constant $K_g$ defined in
(\ref{TST_Def_Kg}).
\end{pf}

The assumption (\ref{TST_Condition_On_Density_Positive_on_the_plane})
ensures that
$K_g>0$, and the regularity constraint
(\ref{TST_Condition_On_Density}) can be verified directly, or using
criteria in
Section~\ref{Section_Understanding_The_Technical_Restriction}.
%
%
\begin{corollary}[(Gaussian zero-mean case)]\label{TST_Corollary_Normal}
If $X\sim\operatorname{MVN}(\boldsymbol{0},\boldsymbol{\Sigma})$,
where $\boldsymbol
{\Sigma}$ is a strictly positive-definite covariance matrix, then
(\ref{TST_First_Order}) holds with
%
%
\begin{equation}
\label{TST_Corollary_Normal_K_g} K_g=C(n_1,n_2,\alpha,\beta)
\frac{\Gamma (\trup{n}{2}
)}{2\uppi^{\trup{n}{2}}|\boldsymbol{\Sigma}|^{1/2}} \int_{-\uppi/2}^{\uppi/2}
\frac{\cos(\omega)^{n-2}} {
(\mathbf{v}(\omega){\boldsymbol\Sigma}^{-1}\mathbf
{v}(\omega)^T  )^{n/2}} \,\mathrm{d}\omega,
\end{equation}
where
$
\mathbf{v}(\omega)
=
\cos(\omega-\omega_0)\mathbbmathbf{I}_{1}
+
\sin(\omega-\omega_0)\mathbbmathbf{I}_{2}$.
\end{corollary}
\begin{pf}
Let $\lambda$ be the smallest eigenvalue of ${\boldsymbol{\Sigma}
}^{-1}$. Note that
$\lambda>0$, which implies that
\[
g(\mathbf{x})\leq\frac{1}{(2\uppi)^{n/2}|{\boldsymbol
\Sigma}|^{1/2}}\mathrm{e}^{-(\trup{\lambda}{2})\|\mathbf{x}\|
^2}<\frac{1}{\|
\mathbf{x}\|^{n+1}}
\]
for $\|\mathbf{x}\|$ large enough. Now, condition
(\ref{TST_Condition_On_Density}) holds according to Lemma~\ref{Fast_Decay_Criterium},
and deriving $K_g$ is a calculus exercise.
\end{pf}

The asymptotic expression for the
distribution tail of the Student two-sample $t$-statistic is obtained by
setting
\[
\alpha=\frac{n_1-1}{n-2} \biggl(\frac{1}{n_1}+\frac{1}{n_2} \biggr)
\quad\mbox{and} \quad \beta=\frac{n_2-1}{n-2} \biggl(\frac{1}{n_1}+
\frac{1}{n_2} \biggr).
\]
For the Gaussian zero-mean case the expression
(\ref{TST_Corollary_Normal_K_g}) then reduces to
%
%
\begin{equation}
\label{TST_Corollary_Normal_K_g_TST} \frac{\Gamma(\trup{n}{2})}{\Gamma(\trup{(n-1)}{2})\sqrt{\uppi
}|\boldsymbol{\Sigma}|^{1/2}} \int_{-\uppi/2}^{\uppi/2}
\frac{\cos(\omega)^{n-2}} {
(\mathbf{v}(\omega){\boldsymbol{\Sigma}}^{-1}\mathbf
{v}(\omega)^T  )^{n/2}} \,\mathrm{d} \omega.
\end{equation}
As expected, if ${\boldsymbol{\Sigma}}=\sigma^2 \mathbf{1}_n$
(recall, $\mathbf{1}_n$ is the identity matrix) and $\sigma^2>0$,
then direct calculation shows that $K_g=1$.
A less trivial case is when the population variances are unequal.
Substituting the diagonal matrix
\[
{{\boldsymbol{\Sigma}}=\operatorname{diag}\bigl\{\underbrace{\sigma
_1^2,\ldots,\sigma _1^2}_{n_1},
\underbrace{\sigma_2^2,\ldots,\sigma_2^2}_{n_2}
\bigr\}}
\]
into (\ref{TST_Corollary_Normal_K_g_TST}), the latter, after some
lengthy algebraic manipulations, takes form
\[
\frac{\Gamma(\trup{n}{2})n_1^{\trup{n}{2}-1}k^{n_2}}{n^{\trup
{(n-1)}{2}}\Gamma(\trup{(n-1)}{2})\sqrt{\uppi}} \biggl[ \int_{-\infty}^{1}
\frac{(1-x)^{n-2}}{(1+ck^2x^2)^{n/2}}\,\mathrm{d}x + \int_1^{\infty}
\frac{(x-1)^{n-2}}{(1+ck^2x^2)^{n/2}}\,\mathrm{d}x \biggr],
\]
where $k={\sigma_1}/{\sigma_2}$ and $c=n_2/n_1$.
The integrals can be computed by resolving the corresponding rational
functions into partial fractions ($n$ is even) or
by expanding brackets in the numerator and integrating by parts ($n$ is
odd). We have computed $K_g$ for sample sizes up to $6$, see
Table~\ref{tab:KgTSTUNEQUALVARIANCES}.
%

%
\begin{table}
\tablewidth=\textwidth
\tabcolsep=0pt
\caption{Constants $K_g$ for the Student two-sample $t$-test,
variances unequal}\label{tab:KgTSTUNEQUALVARIANCES}
\begin{tabular*}{\textwidth}{@{\extracolsep{\fill}}llllll@{}}
\hline
$n_2 \backslash n_1$ & $n_1{=2}$ & $n_1{=3}$ & $n_1{=4}$ & $n_1{=5}$ & $n_1{=6}$
\\
\hline
$n_2{=2}$ & $\frac{k^2+1}{2 k}$ & $\frac{ (2 k^2+3 )^{3/2}}{5
\sqrt{5} k^2}$ & $\frac{ (k^2+2 )^2}{9 k^3}$ & $\frac
{
(2 k^2+5 )^{5/2}}{49 \sqrt{7} k^4}$ & $\frac{ (k^2+3
)^3}{64 k^5}$
\\[5pt]
$n_2{=3}$ & $\frac{ (3 k^2+2 )^{3/2}}{5 \sqrt{5} k}$ &
$\frac
{ (k^2+1 )^2}{4 k^2}$ & $\frac{ (3 k^2+4
)^{5/2}}{49 \sqrt{7} k^3}$ & $\frac{ (3 k^2+5 )^3}{512 k^4}$
& $\frac{ (k^2+2 )^{7/2}}{27 \sqrt{3} k^5}$
\\[5pt]
$n_2{=4}$ & $\frac{ (2 k^2+1 )^2}{9 k}$ & $\frac{ (4
k^2+3 )^{5/2}}{49 \sqrt{7} k^2}$ & $\frac{ (k^2+1
)^3}{8 k^3}$ & $\frac{ (4 k^2+5 )^{7/2}}{2187 k^4}$ &
$\frac
{ (2 k^2+3 )^4}{625 k^5}$
\\[5pt]
$n_2{=5}$ & $\frac{ (5 k^2+2 )^{5/2}}{49 \sqrt{7} k}$ &
$\frac
{ (5 k^2+3 )^3}{512 k^2}$ & $\frac{ (5 k^2+4
)^{7/2}}{2187 k^3}$ & $\frac{ (k^2+1 )^4}{16 k^4}$ & $\frac
{ (5 k^2+6 )^{9/2}}{14\,641 \sqrt{11} k^5}$
\\[5pt]
$n_2{=6}$ & $\frac{ (3 k^2+1 )^3}{64 k}$ & $\frac{  (2
k^2+1 )^{7/2}}{27 \sqrt{3} k^2}$ & $\frac{ (3 k^2+2
)^4}{625 k^3}$ & $\frac{ (6 k^2+5 )^{9/2}}{14\,641 \sqrt{11}
k^4}$ & $\frac{ (k^2+1 )^5}{32 k^5}$
\\
\hline
\end{tabular*}
\end{table}

Note also that for odd sample sizes the exact distribution of the
Student two-sample $t$-statistic
is known, see Ray and Pitman \cite{Ray1961}.

The closed form expressions for (\ref{TST_Corollary_Normal_K_g}) or
(\ref{TST_Corollary_Normal_K_g_TST})
for an arbitrary covariance matrix ${\boldsymbol{\Sigma}}$ is unknown,
but for fixed $n$ one can compute $K_g$ numerically.
In most cases, it is also possible to obtain the exact expression for $K_g$
using {M}athematica \cite{Mathematica2010} software. Examples are
given in
\hyperref[SectionSupplementaryMaterials]{Supplementary Materials}.

\section{Welch statistic}\label{Section_Welch}%

The Welch statistic differs from the Student two-sample $t$-statistic
in that it has $\alpha=1/n_1$ and $\beta=1/n_2$, see the definition
of $T$ in the previous section.
Welch statistic relaxes the assumption of equal variances and its
distribution under the null hypothesis
of equal means is instead approximated by the Student $t$-distribution
with $\nu$ degrees of freedom, where
\[
\nu=\frac{(S_1^2/n_1+S_2^2/n_2)^2}{S_1^4/(n_1^2(n_1-1))+S_2^4/(n_2^2(n_2-1))}
\]
is estimated from the data.
Welch approximation performs poorly in the tail area because it has
wrong asymptotic behavior, cf. Corollary~\ref{TST_Corollary_Normal}.
The accuracy of our asymptotic approximation and its relation to the exact
distribution of the Welch statistic for odd sample sizes, see
Ray and Pitman \cite{Ray1961}, is discussed in
\hyperref[SectionSupplementaryMaterials]{Supplementary Materials}.
We also
study the accuracy of our approximations using simulations, see
Section~\ref{SectionSimulation}.

Finally, Table~\ref{tab:KgWelchTUNEQUALVARIANCES} presents constants
$K_g$ for the Welch statistic under standard assumptions. Here constant
$k$ stands for the ratio $\sigma_1/\sigma_2$.

%
%
\begin{table}
\tablewidth=\textwidth
\tabcolsep=0pt
\caption{Constants $K_g$ for the Welch $t$-test, variances
unequal}\label{tab:KgWelchTUNEQUALVARIANCES}
\begin{tabular*}{\textwidth}{@{\extracolsep{\fill}}lllll@{}}
\hline
$n_2 \backslash n_1$ & $n_1{=2}$ & $n_1{=3}$ & $n_1{=4}$ & $n_1{=5}$
\\
\hline
$n_2{=2}$ & $\frac{k^2+1}{2 k}$ & $\frac{ (2 k^2+3 )^{3/2}}{9
k^2}$ & $\frac{3 \sqrt{(\trup{3}{2})}  (k^2+2 )^2}{16 k^3}$ &
$\frac{4  (2 k^2+5 )^{5/2}}{125 k^4}$
\\[5pt]
$n_2{=3}$ & $\frac{ (3 k^2+2 )^{3/2}}{9 k}$ & $\frac{
(k^2+1 )^2}{4 k^2}$ & $\frac{ (3 k^2+4 )^{5/2}}{50
\sqrt{5} k^3}$ & $\frac{4  (3 k^2+5 )^3}{1215 k^4}$
\\[5pt]
$n_2{=4}$ & $\frac{3 \sqrt{(\trup{3}{2})}  (2 k^2+1 )^2}{16 k}$
& $\frac{ (4 k^2+3 )^{5/2}}{50 \sqrt{5} k^2}$ & $\frac
{
(k^2+1 )^3}{8 k^3}$ & $\frac{3 \sqrt{(\trup{3}{35})}  (4
k^2+5 )^{7/2}}{1715 k^4}$
\\[5pt]
$n_2{=5}$ & $\frac{4  (5 k^2+2 )^{5/2}}{125 k}$ & $\frac{4
(5 k^2+3 )^3}{1215 k^2}$ & $\frac{3 \sqrt{(\trup{3}{35})}
(5 k^2+4 )^{7/2}}{1715 k^3}$ & $\frac{ (k^2+1
)^4}{16 k^4}$
\\[5pt]
$n_2{=6}$ & $\frac{25 \sqrt{(\trup{5}{3})}  (3 k^2+1 )^3}{216
k}$ & $\frac{25 \sqrt{(\trup{5}{7})}  (2 k^2+1 )^{7/2}}{343
k^2}$ & $\frac{25 \sqrt{(\trup{5}{2})}  (3 k^2+2 )^4}{16\,384
k^3}$ & $\frac{4  (6 k^2+5 )^{9/2}}{177\,147 k^4}$
\\
\hline
\end{tabular*}
\end{table}

%
\section{$F$-statistic}\label{Section_F}%

In this section, we study the tails of the distribution of an
$F$-statistic for testing the equality of variances. Similar results
can also be obtained for an $F$-test used in one-way ANOVA, lack-of-fit
sum of squares, and when comparing two nested linear models in
regression analysis.
Define random vectors $\mathbf{X}=(X_1,X_2,\ldots,X_{n_1})$ and
$\mathbf{Y}=(Y_1,Y_2,\ldots,Y_{n_2})$, $n_1\geq2$ and $n_2\geq2$,
and let $g(\mathbf{x},\mathbf{y})$ be the joint density of the vector
$(\mathbf{X},\mathbf{Y})$. Now set $n=n_1+n_2$ and define
\[
T = S_1^2/S_2^2,
\]
where $S_1$ and $S_2$ are the sample variances of $\mathbf{X}$ and
$\mathbf{Y}$, respectively.
Let $s_1(\mathbf{x})$ denote the sample standard deviation of the
vector $\mathbf{x}\in\mathbb{R}^{n_1}$ and define the unit vector
$\mathbbmathbf{I}= (1/{\sqrt{n_2}},1/{\sqrt{n_2}},\ldots
,\allowbreak  1/{\sqrt{n_2}} )$.
We assume that
%
%
\begin{equation}
\label{F_Condition_On_Density_Positive_on_the_plane} s_1(\mathbf{x}) g ( \mathbf{x},r\mathbbmathbf{I} )>0
\end{equation}
for some $\mathbf{x}$ and $r$,
and that the integral
%
%
\begin{equation}
\label{F_Condition_On_Density} \mathop{\idotsint}_{\mathbb{R}^{n_1}} s_1(
\mathbf{x})^{n_2-1} \int_{-\infty}^{\infty} \mathop{
\max_{\|{\boldsymbol{\xi}}\|<\varepsilon,}}_{
{\boldsymbol{\xi}}\in L^\perp}g \bigl( \mathbf{x},r\mathbbmathbf
{I}+s_1(\mathbf{x}) {\boldsymbol\xi} \bigr)\,\mathrm{d}r \,
\mathrm{d}\mathbf{x}
\end{equation}
is finite for some $\varepsilon>0$, where $L$ is a linear subspace
spanned by vector $\mathbbmathbf{I}$ and $L^\perp$ is its
orthogonal complement. Finally, define the constant
%
%
\begin{equation}
\label{F_Def_Kg} K_g= \frac{\Gamma (\trup{(n_1-1)}{2} ) (\uppi
(n_1-1)
)^{\trup{(n_2-1)}{2}}}{\Gamma (\trup{(n-2)}{2} )} \mathop{
\idotsint}_{\mathbb{R}^{n_1}} s_1(\mathbf{x})^{n_2-1} \int
_{-\infty}^{\infty} g (\mathbf{x},r\mathbbmathbf{I} )\,
\mathrm{d}r \,\mathrm {d}\mathbf{x}.
\end{equation}

%
\begin{theorem}
\label{FFirstOrderTheorem}
If $g$ is continuous and satisfies
(\ref{F_Condition_On_Density_Positive_on_the_plane}) and
(\ref{F_Condition_On_Density}), then
%
%
\begin{equation}
\label{F_First_Order} \frac{\mathbf{P}({T>u})}{F_{n_1-1,n_2-1}(u)}=K_g+\mathrm{o}(1) \qquad\mbox{as
} u\to\infty,
\end{equation}
where $F_{n_1-1,n_2-1}(u)$ is the tail of the $F$-distribution with
parameters $n_1-1$ and $n_2-1$ and $0<K_g=K(g)<\infty$.
\end{theorem}
%
%
\begin{corollary}[(Gaussian zero-mean case, independent samples)]\label
{F_First_Order_Gaussian_Corollary}
If $X$ and $Y$ are independent zero-mean Gaussian random vectors with
strictly non-degenerate covariance matrices
${\boldsymbol{\Sigma}_1}$ and ${\boldsymbol{\Sigma}_2}$, then
(\ref{F_First_Order}) holds with
%
%
\begin{equation}
\label{F_First_Order_Gaussian} K_g = C \mathop{\idotsint}_{\mathbb{R}^{n_1}}
\frac{s_1(\mathbf
{x})^{n_2-1}}{
(1+\mathbf{x}{\boldsymbol\Sigma}^{-1}_1\mathbf{x}^T )^{n/2}} \,\mathrm{d}\mathbf{x},
\end{equation}
where the constant $C$ is given by
\[
C=\frac{(n-2)(n_1-1)^{\trup{(n_2-1)}{2}}\Gamma (\trup
{(n_1-1)}{2} )|\mathbbmathbf{I}{\boldsymbol\Sigma
}_2\mathbbmathbf{I}^T|^{1/2}} {
2\uppi^{\trup{(n_1+1)}{2}} |{\boldsymbol\Sigma}_1|^{1/2}
|{\boldsymbol
\Sigma}_2|^{1/2}}.
\]
\end{corollary}

The proofs of Theorem~\ref{FFirstOrderTheorem} and Corollary~\ref{F_First_Order_Gaussian_Corollary} are given in
\hyperref[SectionSupplementaryMaterials]{Supplementary Materials}.
Now consider the asymptotic power of the $F$-statistic.

%
\begin{corollary}[(Asymptotic power)]\label{F_First_Order_Power_Corollary}
If $X$ and $Y$ are independent zero-mean Gaussian random vectors with
covariance matrices $\sigma_1^2 \mathbf{1}_{n_1}$ and $\sigma_2^2
\mathbf{1}_{n_2}$, $\sigma_1^2+\sigma_2^2>0$, then
%
%
\begin{equation}
\label{F_First_Order_Power} \lim_{u\to\infty}\frac{\mathbf
{P}({T>u})}{F_{n_1-1,n_2-1}(u)}= \biggl(
\frac
{\sigma_1}{\sigma_2} \biggr)^{n_2-1}.
\end{equation}
\end{corollary}
\begin{pf}
Changing variables $\mathbf{x}=\sigma_1 B\mathbf{y}$, where $B$ is
an orthogonal operator such that $B\mathbf{e}_{\mathbf{n}_\mathbf
{1}}= (1/{\sqrt {n_1}},1/{\sqrt{n_1}},\ldots,1/{\sqrt{n_1}} )$, the integral
on the
right-hand side of (\ref{F_First_Order_Gaussian}) takes form
\[
\sigma_1^{n-1} \biggl(\frac{1}{n_1-1}
\biggr)^{\trup
{(n_2-1)}{2}}\mathop{\idotsint}_{\mathbb{R}^{n_1}}\frac{ (\|
\mathbf{y}\|
^2-y_{n_1}^2 )^{\trup{(n_2-1)}{2}}}{ (1+\|\mathbf{y}\|
^2 )^{n/2}}\,
\mathrm{d}\mathbf{y},
\]
and is then evaluated by passing to spherical coordinates.
\end{pf}

A careful reader may note that (\ref{F_First_Order_Power}) follows
from the asymptotic expansion of
$
\mathbf{P}({T>u})=F_{n_1-1,n_2-1} ((\sigma_2/\sigma_1)^2
u )
$
in terms of
$F_{n_1-1,n_2-1} (u )$. Our aim was just to show that
despite the seeming complexity of the expression (\ref{F_Def_Kg}), the
constant $K_g$
can be evaluated directly, at least for some standard densities.
It is also possible to compute $K_g$ numerically, see the
{MATLAB} \cite{Matlab2010} scripts in
\hyperref[SectionSupplementaryMaterials]{Supplementary Materials}.

\section{Second and higher order approximations}%
\label{Section_Second_Order_Approximations}%

In this section, we discuss the speed of convergence in
Theorem~\ref{IntroductionMainTheorem}.
Let $T$ be one of the test statistics defined in
Sections~\ref{Section_One_Sample}, \ref{Section_Two_Sample} and~\ref
{Section_F}
and let $t_k(u)$ be the Student $t$-distribution tail with $k$ degrees
of freedom
and $F_{m,k}(u)$ be the $F$-distribution tail with
parameters $m$ and $k$.
For an arbitrary continuous multivariate density $g=g_1(\mathbf{x})$,
assume that conditions
(\ref{OST_Condition_On_Density}), (\ref{TST_Condition_On_Density})
and (\ref{F_Condition_On_Density})
hold, and define the constant $K_g$ by (\ref{OST_Def_Kg}),
(\ref{TST_Def_Kg}) and (\ref{F_Def_Kg}) for the three tests respectively.
For the Student $t$-statistic the function $G(\mathbf{t})$ is given
by (\ref{OST_Iterated_int}) and (\ref{TST_Iterated_int}), and for the
$F$-statistic see the corresponding formula in the proof of
Theorem~\ref{FFirstOrderTheorem} in
\hyperref[SectionSupplementaryMaterials]{Supplementary Materials}.
Finally, with the standard notation $\nabla f$ for the gradient of a
scalar function $f$, and a parameter $\alpha$ which can take values
$1$ or $2$,
define
%
%
\begin{equation}
d_{\alpha,m,k}(u)=\frac{1}{u^{\trup{\alpha(k+1)}{2}}} \biggl[ C_1\sup
_{{\|\mathbf{x}\|}\leq u^{-\trup{\alpha}{2}}}\bigl\llVert \nabla G(\mathbf{x}) \bigr\rrVert +
C_2\frac{K_{g}}{\alpha}\frac{1}{u^{\trup{\alpha}{2}}} \biggr],
\end{equation}
where the constants $C_1$, $C_2$ (which depend on $m$ and $k$) are
given in Lemma~\ref{IOFOSB_Base_Lemma}(B).
%
%
\begin{lemma}[(Absolute error bound)]\label{AERROR_Bound}
If $G(\mathbf{t})$ is differentiable in some neighborhood of zero,
then for any $u>0$
the following inequalities
\begin{eqnarray*}
\bigl\llvert \mathbf{P}({T>u})-K_g t_{n-1}(u)\bigr\rrvert
&\leq& d_{2,1,n-1}(u),
\\
\bigl\llvert \mathbf{P}({T>u})-K_g t_{n-2}(u)\bigr\rrvert
&\leq& d_{2,1,n-2}(u),
\\
\bigl\llvert \mathbf{P}({T>u})-K_g F_{n_1-1,n_2-1}(u)\bigr\rrvert
&\leq& d_{1,n_1-1,n_2-1}(u),
\end{eqnarray*}
hold for the Student one- and two-sample $t$- and $F$-statistics accordingly.
\end{lemma}

\begin{pf}
The first two inequalities follow from (\ref{OST_First_Order}),
(\ref{TST_First_Order}) and Corollary~\ref{IOFOSB_Base_Lemma_Corollary},
and
for the $F$-statistic we use
Lemma~\ref{IOFOSB_Base_Lemma}(B)
with $\alpha=1$ and
$\sqrt{u}$ instead of $u$.
\end{pf}

Below follows the asymptotic formula for the relative error.
For convenience, we denote
the distribution tail of $T$ under the null hypothesis $H_0\dvtx
g_0\sim
\operatorname{MVN}(\mathbf{0},\sigma^2\mathbf{1}_n)$ by $t(u)$.
%
%
\begin{lemma}[(Relative error decrease rate)]\label{RERROR_Bound}
If $G(\mathbf{t})$ is twice differentiable in some neighborhood of
zero, then
\[
\frac{\mathbf{P}({T>u})-K_g t(u)}{\mathbf{P}({T>u})}=\frac
{C_3}{u^\alpha}\bigl(1+\mathrm{o} (1 )\bigr),
\]
where
\[
C_3=\frac{\alpha k B (\trup{m}{2},\trup{k}{2} )}{2
(\trup{k}{m} )^{k/2}} \frac{L_{G,\alpha}}{K_g},
\]
the triple $(\alpha,m,k)$ is set to
$(2,1,n-1)$, $(2,1,n-2)$ and $(1,n_1,n_2)$ for the Student one- and
two-sample $t$- and $F$-statistics, respectively, and
the constant $L_{G,\alpha}$ is defined in
Lemma~\textup{\ref{IOFOSB_Base_Lemma}(C)}.
\end{lemma}
\begin{pf}
The result follows from formulas (\ref{OST_First_Order}),
(\ref{TST_First_Order}) and (\ref{F_First_Order})
for $\mathbf{P}({T>u})$, Lemma~\ref{IOFOSB_Base_Lemma}(C)
and formula (\ref{IOFOSB_t_distribution_tail_asymptotics}).
\end{pf}

The bounds and asymptotic expressions for
the case of an arbitrary null hypothesis $H_0$ are derived using basic calculus:
\begin{eqnarray*}
&&\mathbf{P}({T>u|H_1})-({K_{g_1}}/{K_{g_0}})
\times\mathbf{P}({T>u|H_0})
\\
&&\quad= \bigl(\mathbf{P}({T>u|H_1})-K_{g_1}t(u) \bigr) -
({K_{g_1}}/{K_{g_0}})\times \bigl( \mathbf{P}({T>u|H_0})-K_{g_0}t(u)
\bigr),
\end{eqnarray*}
and the absolute error of the approximation
(\ref{Introduction_Main_result}) is thus bounded by the linear combination
of the absolute errors considered in Lemma~\ref{AERROR_Bound} above.

For the relative error, we replace the two probabilities $\mathbf
{P}({T>u|H_1})$ and $\mathbf{P}({T>u|H_0})$
by their second order expansions given by Lemma~\ref{IOFOSB_Base_Lemma}(C),
and then use (\ref{IOFOSB_t_distribution_tail_asymptotics}).
Lemma~\ref{IOFOSB_Base_Lemma}
can also be generalized to obtain higher order series expansion for
$\mathbf{P}({T>u})$ as $u\to\infty$.

\section{Simulation study}\label{SectionSimulation}%

Let $T$ be one of the test statistics considered in the previous
sections and
$t(u)$ be the distribution tail of $T$ under $H_0\dvt g\sim
\operatorname{MVN}(\mathbf
{0},\mathbf{1}_n)$.
Next, we choose the sample size, specify the density $g(\mathbf{x})$, and
simulate $N$ random vectors $\mathbf{X}\sim g$. For each vector
$\mathbf{X}$, we compute $t^*=T(\mathbf{X})$, the value of the test
statistic $T$,
and two $p$-values $p^{R}=t(t^*)$ and $p^{C}=K_g t(t^*)=K_gp^R$. Finally,
we plot the empirical CDF of $p^{R}$ and $p^{C}$ over the range $I(r)=[0,1/r]$,
where the \textit{Zoom Factor} (Z.F.) parameter $r$ determines the
tail region of interest. Here $N=10\,000\times r$ so that
$I(r)$ contains approximately $10\,000$ $p$-values (as if they were
uniformly distributed) -- this is to ensure that the tails of the
distribution of the $p$-values $p^R$ and $p^C$
are equally well approximated by the corresponding CDFs in all the tail
regions. The letters ``R'' and ``C'' in the notation for the $p$-values
stand for ``Raw'', that is, computed using $t(u)$, and ``Corrected'',
that is, computed using $K_gt(u)$.\vadjust{\goodbreak}

For the i.i.d. case, let $h(x)$, the marginal density of the vector
$\mathbf{X}$, be either $\operatorname{Uniform}(-1,1)$, \textit{Standard
normal}, \textit{Centered exponential}, \textit{Cauchy}, or
$t$-density with $2$ or $5$ degrees of freedom.
The constant $K_g$ was either evaluated explicitly in {M}athematica
\cite{Mathematica2010}
or computed numerically in
{MATLAB} \cite{Matlab2010},
see \hyperref[SectionSupplementaryMaterials]{Supplementary Materials}.
Figures~\ref{fig:SimulationOST}, \ref{fig:SimulationTST} and~\ref{fig:SimulationF} in \hyperref
[SectionFiguresandTables]{Appendix B}
show empirical CDFs for
different sample sizes and Zoom Factor $r$ varying between $20$ and $1\,
000\,000$.
One can see that our approximations are very accurate in the tail
regions for all the three test statistics, all sample sizes, and
densities $h(x)$ considered in the study.
Note also that the convergence speed is better for smaller sample sizes --
this is in agreement with the bounds for the absolute error in Lemma~\ref{AERROR_Bound}, see Section~\ref{Section_Second_Order_Approximations}.

Next, we computed $p$-values for the Welch statistic and compared them
with the $p$-values obtained using the expression
(\ref{TST_Corollary_Normal_K_g}) in Corollary~\ref{TST_Corollary_Normal}.
Here ``Raw'' $p$-values are obtained using the Welch approximation and
the notation is $p^{W}$.
According to the plots in the top row of Figure~\ref
{fig:SimulationWELCH}, it may seem that the $p$-values $p^{W}$ are
uniformly distributed.
However, if one ``zooms in'' to the tail region, see the plots in the
middle row of Figure~\ref{fig:SimulationWELCH}, it is clear that the
$p$-values obtained using Welch approximation deviate significantly from
the theoretical uniform distribution, while the corrected $p$-values
$p^C$ follow the diagonal line precisely.
The advantage of using our tail approximations is fully convincing at
Zoom Factor $100\,000$, see the bottom row of Figure~\ref{fig:SimulationWELCH}.

Finally, we made similar plots for even more peculiar scenarios where
the data was dependent and non-stationary, see, for example,
Figure~\ref{fig:SimulationMVN}. Our approximations were very accurate
in all
considered cases.

%
\begin{appendix}

\section{Supplementary theorems and lemmas}\label
{Section_FurtherRemarks}%

This Appendix is split into two parts -- the first one introduces the
key lemma which is used in Sections~\ref{Section_One_Sample},
\ref{Section_Two_Sample} and~\ref{Section_F}, and the
second contains useful notes on the regularity constraints (replacing
them by simpler criteria that can be used in practice) and shows how to
weaken the assumption of continuity of the density $g(\mathbf{x})$.

\subsection{Asymptotic behavior of an integral of a continuous
function over a shrinking ball}\label{Section_Shrinking_Balls}%

It was shown that the tails of the distribution of the Student one- and
two-sample $t$-, Welch, and $F$-statistics are determined by the
asymptotic behavior of an integral of some function (different for each
of the tests) over a shrinking ball.

Let $G(\mathbf{x})$, $\mathbf{x}\in\mathbb{R}^k$ be some
real-valued function
and consider the asymptotic behavior of
\setcounter{equation}{22}
\renewcommand{\theequation}{\arabic{equation}}
%
\begin{equation}
\label{IOFOSB_The_Integral} F(u)=\mathop{\idotsint}_{\sum x_i^2<u^{-2}} G(\mathbf{x}) \,
\mathrm{d}\mathbf{x}
\end{equation}
for fixed $k$ and $u\to\infty$.\vadjust{\goodbreak}
%
%
\begin{lemma}
\label{IOFOSB_Base_Lemma}
Set $f(u)=\alpha^{-1}F_{m,k}(u^2)$, where $F_{m,k}(\cdot)$ is the
tail of the $F$-distribution with $m\geq1$ and $k\geq2$ degrees of freedom,
and let $\operatorname{Vol}(B_k)$ be the volume of the unit $k$-ball $B_k$
and $B(x,y)$ be the Beta function. The parameters $\alpha$ and $m$
will be set later.
With the above notation we have:
\begin{enumerate}[(B)]
\item[(A)]If $G$ is continuous at zero, then
%
%
\begin{equation}
\label{IOFOSB_First_Order} \frac{F(u)}{f(u)}=K_{G,\alpha} +\mathrm{o}(1) \qquad
\mbox{as } u\to\infty,
\end{equation}
where
%
%
\begin{equation}
\label{IOFOSB_Def_K_G} K_{G,\alpha}= \frac{\alpha k B (\trup{m}{2},\trup{k}{2} )}{2
(\trup{k}{m} )^{k/2}} \operatorname{Vol}(B_k)
G(\mathbf{0}).
\end{equation}
\item[(B)]If $G$ is differentiable in some neighborhood of zero,
then for any $u>0$
%
%
\begin{equation}
\label{IOFOSB_Absolute_Error_Bound} \bigl\llvert F(u)-K_{G,\alpha} f(u)\bigr\rrvert \leq
\frac{C_1}{u^{k+1}} \sup_{{\|\mathbf{x}\|}\leq u^{-1}}\bigl\llVert \nabla G(\mathbf{x})
\bigr\rrVert + C_2\frac{K_{G,\alpha}} {
\alpha}\frac{1}{u^{k+2}},
\end{equation}
where
%
%
\begin{equation}
\label{IOFOSB_Absolute_Error_Bound_Def_C} C_1=\operatorname{Vol}(B_k) \quad
\mbox{and} \quad C_2=\frac{k(k+m)}{m(k+2)}\frac{ (\trup{k}{m}
)^{k/2}}{B (\trup{m}{2},\trup{k}{2} )},
\end{equation}
and $\nabla G(\mathbf{x})$ is a gradient of $G$ evaluated at point
$\mathbf{x}$.
\item[(C)]If $G$ is twice differentiable in some neighborhood
of zero, then
%
%
\begin{equation}
\label{IOFOSB_Second_Order} u^{k+2} \bigl(F(u)-K_{G,\alpha} f(u)
\bigr)=L_{G,\alpha}+\mathrm{o}(1) \qquad\mbox{as } u\to\infty,
\end{equation}
where
\[
L_{G,\alpha}=C_1 \frac{\operatorname{tr} (\operatorname{Hess} (G(\mathbf
{0}) ) )}{2(k+2)} - C_2
\frac{K_{G,\alpha}}{\alpha} ,
\]
$\operatorname{tr}(A)$ is the trace of a square matrix $A$, and
$\operatorname{Hess} (G(\mathbf
{x}) )$
is the Hessian matrix of $G$ evaluated at point $\mathbf{x}$.
Constants $C_1$ and $C_2$ are given by (\ref{IOFOSB_Absolute_Error_Bound_Def_C}).
\end{enumerate}
\end{lemma}
\begin{pf}
The first statement follows from the asymptotic expansion for the
$F$-distribution tail
%
\begin{equation}
\label{IOFOSB_t_distribution_tail_asymptotics} f(u)= \frac{2 (\trup{k}{m} )^{k/2}}{\alpha k B (\trup
{m}{2},\trup{k}{2} )} \biggl[\frac{1}{u^{k}} -
\frac{k^2(k+m)}{2m(k+2)}\frac{1}{u^{k+2}} \biggr] + \mathrm{o} \biggl(
\frac{1}{u^{k+2}} \biggr).
\end{equation}
Indeed, changing variables $\mathbf{x}=\mathbf{y}/u$ we write
%
%
\begin{equation}
\label{IOFOSB_The_Integral_After_Variable_Change} F(u)=\mathop{\idotsint}_{\sum x_i^2<u^{-2}} G(\mathbf{x}) \,\mathrm
{d}\mathbf{x}= \frac{1}{u^{k}}\mathop{\idotsint}_{B_k} G(
\mathbf{y}/u) \,\mathrm {d}\mathbf{y}.
\end{equation}
Continuity of $G$ at zero implies uniform convergence of
$G(\mathbf{y}/u)$ to $G(\mathbf{0})$ over the ball ${B_k}$,
and thus
%
%
\begin{equation}
\label{IOFOSB_natural_first_order_approximation} F(u)=\operatorname{Vol} (B_k )G(\mathbf{0})
\frac{1}{u^{k}} \bigl(1+\mathrm{o}(1) \bigr).
\end{equation}
Dividing (\ref{IOFOSB_natural_first_order_approximation}) by
(\ref{IOFOSB_t_distribution_tail_asymptotics})
we get that the value of $K_{G,\alpha}$ in
(\ref{IOFOSB_First_Order}) coincides with (\ref{IOFOSB_Def_K_G}).

Now assume $G$ is differentiable in some neighborhood of zero and consider
the Lagrange form of the Taylor expansion of $G(\mathbf{y}/u)$.
The latter and (\ref{IOFOSB_The_Integral_After_Variable_Change}) give
\begin{eqnarray*}
\bigl\llvert F(u)-K_{G,\alpha} f(u)\bigr\rrvert &\leq& \frac{1}{u^k}
\bigl\llvert \operatorname{Vol} (B_k )G(\mathbf{0})-u^{k}K_{G,\alpha}
f(u)\bigr\rrvert
\\
&&{}+\frac{1}{u^{k+1}}\biggl\llvert \mathop{\idotsint}_{B_k} \nabla G
\bigl(\xi (\mathbf{y})\mathbf{y} \bigr)\mathbf{y}^{T} \,\mathrm{d}
\mathbf {y}\biggr\rrvert ,
\end{eqnarray*}
where $0\leq\xi(\mathbf{y})\leq1/u$.
The second summand in the right-hand side of the above inequality is
bounded by
\[
\frac{1}{u^{k+1}}\operatorname{Vol} (B_k )\sup
_{B_k}\bigl \| \nabla G(\mathbf{x}/u)\bigr \|,
\]
and the bound for the remaining summand follows from
(\ref{IOFOSB_t_distribution_tail_asymptotics}), where
we note that $f(u)$ is bounded by the two successive partial sums in
its alternated series (\ref{IOFOSB_t_distribution_tail_asymptotics})
and that the factors before $\operatorname{Vol} (B_k
)G(\mathbf
{0})$ in the expression for $K_{G,\alpha}$ and before the square
brackets in (\ref{IOFOSB_t_distribution_tail_asymptotics}) cancel out.
The last step is to use formulas (\ref{IOFOSB_Def_K_G}) and
(\ref{IOFOSB_Absolute_Error_Bound_Def_C}) to express $\operatorname
{Vol}
(B_k )G(\mathbf{0})$ in terms of $K_{G,\alpha}$ and~$C_2$.

We move on to the proof of (\ref{IOFOSB_Second_Order}). Taylor
expansion for $G(\mathbf{y}/u)$ yields
\begin{eqnarray*}
F(u) &=& \frac{1}{u^k}\operatorname{Vol} (B_k )G(\mathbf{0})
\\
&&{}+ \frac{1}{u^{k+2}}\mathop{\idotsint}_{B_k} \frac{\mathbf
{y}\operatorname{Hess} (G
(\mathbf{0} ) )\mathbf{y}^{T}}{2}
\,\mathrm{d}\mathbf{y} + \mathrm{o} \biggl(\frac{1}{u^{k+2}} \biggr),
\end{eqnarray*}
where we took into account that
the integral of the odd function $\nabla G(\mathbf{0})\mathbf{y}$
over the ball $B_k$ is zero.
Neglecting odd terms in $\mathbf{y}\operatorname{Hess} (G
(\mathbf
{0} ) )\mathbf{y}^{T}$ we have
\begin{eqnarray*}
\mathop{\idotsint}_{B_k} \mathbf{y}\operatorname{Hess} \bigl(G (
\mathbf{0} ) \bigr)\mathbf{y}^{T}\,\mathrm{d}\mathbf{y}&=& \sum
\mathop{\idotsint}_{B_k} \frac{\partial^2G(\mathbf
{0})}{\partial
^2y_i}y_i^2
\,\mathrm{d}\mathbf{y}
\\
&=& \biggl(\sum\frac{\partial^2G(\mathbf{0})}{\partial^2y_i} \biggr)\mathop{\idotsint}_{B_k}
\frac{\sum y_i^2}{k} \,\mathrm{d}\mathbf{y}
\\
&=& \operatorname{Vol} (B_k
)\frac{\operatorname
{tr}(\operatorname{Hess}(G(\mathbf{0}))}{k+2},
\end{eqnarray*}
where the last integral was computed using spherical coordinates.
Substituting the second order Taylor expansion for $F(u)$ and
expression for $f(u)$ in (\ref{IOFOSB_t_distribution_tail_asymptotics})
into
the left-hand side of (\ref{IOFOSB_Second_Order}) we get the constant
$L_{G,\alpha}$.
\end{pf}

Note that the expression ${\alpha}^{-1}{K_{G,\alpha}}$ does not
depend on $\alpha$ and thus the right-hand side of
(\ref{IOFOSB_Absolute_Error_Bound}) and (\ref{IOFOSB_Second_Order})
depends
only on the integrand $G$ in (\ref{IOFOSB_The_Integral}) and
parameters $m$ and $k$.\vspace*{-2pt}
%
%
\begin{corollary}
\label{IOFOSB_Base_Lemma_Corollary}
Let $t_k(u)$ be the Student $t$-distribution tail with $k$ degrees of
freedom. If $G$ is continuous at zero, then
\[
\frac{F(u)}{t_k(u)}=K_{G,2} +\mathrm{o}(1) \qquad\mbox{as } u\to
\infty,
\]
where $K_{G,2}$ is given by \textup{(\ref{IOFOSB_Def_K_G})} with $m=1$.
Statements \textup{(B)} and \textup{(C)} also hold for $f(u)=t_k(u)$,
provided $m=1$ and
$\alpha=2$.\vspace*{-2pt}
\end{corollary}
\begin{pf}
Note that $t_k(u)=\frac{1}{2}F_{1,k}(u^2)$ and apply
Lemma~\ref{IOFOSB_Base_Lemma}.\vspace*{-2pt}
\end{pf}
%
\subsection{A note on the regularity constraints and the
continuity assumption}%
\label{Section_Understanding_The_Technical_Restriction}%

The aim of this section is to replace the technical constraints
(\ref{OST_Condition_On_Density}), (\ref{TST_Condition_On_Density})
and (\ref{F_Condition_On_Density})
of Theorems~\ref{OST_First_Order_Theorem},
\ref{TST_First_Order_Theorem} and~\ref{FFirstOrderTheorem}
by simpler criteria, and to weaken the assumption of continuity of the
multivariate density $g(\mathbf{x})$ of the data vector $\mathbf{X}$.

The nature of the regularity constraints
(\ref{OST_Condition_On_Density}), (\ref{TST_Condition_On_Density}) and
(\ref{F_Condition_On_Density})
becomes clear if one notes that all the proofs
share a common part, which is to apply
Lemma~\ref{IOFOSB_Base_Lemma}(A) or Corollary~\ref
{IOFOSB_Base_Lemma_Corollary} to the
representation for the distribution tail of the test statistic $T$, see
(\ref{OST_Iterated_int}) and (\ref{TST_Iterated_int}), and then to
use dominated convergence theorem to show that
the corresponding function $G(\mathbf{t})$ is continuous at zero.
The only purpose of the regularity constraints is to ensure that the
limiting and integration operations are interchangeable,
and that the resulting constant $K_g$ is finite.
Omitting the regularity assumptions (\ref{OST_Condition_On_Density}),
(\ref{TST_Condition_On_Density}) and (\ref{F_Condition_On_Density})
we immediately obtain\vspace*{-2pt}

%
\begin{theorem}[(``liminf'' analogue of Theorems~\ref
{OST_First_Order_Theorem}, \ref{TST_First_Order_Theorem} and~\ref
{FFirstOrderTheorem})]\label{Relaxing_Assumptions_LimInf_Theorem}
Let $T$ be the Student one- or two-sample $t$-statistic or an
$F$-statistic and let $t(u)$ be the distribution tail of $T$ under the
null hypothesis $H_0\dvt g\sim\operatorname{MVN}(\mathbf{0},\sigma
^2\mathbf{1}_n)$,
where $\sigma^2>0$ and $\mathbf{1}_n$ is the identity matrix.
If $g$ is continuous, then
\[
\liminf_{u\to\infty}\frac{\mathbf{P}({T>u})}{t(u)}\geq K_g,
\]
where the constant $K_g$ is given by (\ref{OST_Def_Kg}),
(\ref{TST_Def_Kg}) and (\ref{F_Def_Kg}) accordingly,
though it may not necessarily be finite.\vspace*{-2pt}
\end{theorem}

Next, we give the sufficient (but not necessary) conditions for the
regularity constraints of Theorems~\ref{OST_First_Order_Theorem},
\ref{TST_First_Order_Theorem} and~\ref{FFirstOrderTheorem} to
hold. One
may expect that formulas (\ref{OST_First_Order}),
(\ref{TST_First_Order}) and (\ref{F_First_Order}) hold simply when
$g$ is
continuous and $K_g$ is finite, but proving or disproving
this claim is not easy and it remains an open problem.\vadjust{\goodbreak}

%
\begin{lemma}
\label{Fast_Decay_Criterium}
If $g(\mathbf{x})$ is bounded and there exist positive constants $R$,
$C$ and $\delta$ such that
%
%
\begin{equation}
\label{Fast_Decay_Criterium_Bound} g(\mathbf{x})\leq\frac{C}{\|\mathbf{x}\|^{n+\delta}} \qquad\mbox {for } \|
\mathbf{x}\|>R,
\end{equation}
then the assumptions (\ref{OST_Condition_On_Density}),
(\ref{TST_Condition_On_Density}) and (\ref{F_Condition_On_Density})
of Theorems~\ref{OST_First_Order_Theorem}, \ref{TST_First_Order_Theorem}
and~\ref{FFirstOrderTheorem} hold.
\end{lemma}
\begin{pf}
The integrals in (\ref{OST_Condition_On_Density}),
(\ref{TST_Condition_On_Density}) and (\ref{F_Condition_On_Density})
will be
estimated by partitioning the integration domain into several disjoint
parts $D_i$ and $D^{*}_j$ and analyzing the integrals over these sets
separately.
For non-compact domains $D^{*}_j$ the integrand will be estimated from
above using the bound (\ref{Fast_Decay_Criterium_Bound}) and showing
that this bound is integrable.
The integrability over the compact domains $D_i$ follows from the fact
that $g(\mathbf{x})$ is bounded.
In the notation below let $G(r)$, $G(\omega,r)$ and $G(\mathbf{x},r)$
be the integrands in
(\ref{OST_Condition_On_Density}), (\ref{TST_Condition_On_Density})
and (\ref{F_Condition_On_Density}) accordingly.

\noindent
\textit{Student's one-sample $t$-statistic}: Set $D_1=[0,R]$ and
$D^{*}_1=[R,\infty]$.
Since $\mathbbmathbf{I}$ and ${\boldsymbol{\xi}}$ are orthogonal
and taking into account that $\|\mathbbmathbf{I}\|=1$ we have
$\|r
(
\mathbbmathbf{I}
+
{\boldsymbol\xi}
)\|^2
=r^2(1+\|{\boldsymbol\xi}\|^2)\geq r^2$,
and the bound (\ref{Fast_Decay_Criterium_Bound}) gives
\[
\int_{D_1^{*}} G(r) \,\mathrm{d}r < \int_{R}^{\infty}
\frac{C}{r^{1+\delta}} \,\mathrm{d}r<\infty.
\]

\noindent
\textit{Student's two-sample $t$-statistic}:
Setting
$
D_1=[-\uppi/2,\uppi/2]\times[0,R]$ and $D^{*}_1=[-\uppi/2,\uppi
/2]\times
[R,\infty]
$
and noting that $\mathbbmathbf{I}_{1}$, $\mathbbmathbf{I}_{2}$ and
${\boldsymbol\xi}$ are mutually orthogonal we get
\[
\bigl \|r \bigl( \cos(\omega-\omega_0)\mathbbmathbf{I}_{1} +
\sin(\omega-\omega_0)\mathbbmathbf{I}_{2} + {\boldsymbol
\xi} \bigr)\bigr \|^2=r^2\bigl(1+\|{\boldsymbol\xi}
\|^2\bigr)\geq r^2,
\]
where we used the fact that $\|\mathbbmathbf{I}_{1}\|=\|\mathbbmathbf
{I}_{2}\|=1$. Now
the bound (\ref{Fast_Decay_Criterium_Bound}) implies
\[
\int_{D_1^{*}} G(\omega,r) \,\mathrm{d}r < \int
_{-\uppi/2}^{\uppi/2} \cos(\omega)^{n-2} \,\mathrm{d}
\omega\times \int_{R}^{\infty} \frac{C}{r^{1+\delta}}<
\infty.
\]

\noindent
\textit{$F$-statistic}: Consider the following partition of $\mathbb
{R}^{n_1+1}\dvtx
D_1= \{ (\mathbf{x},r)\dvtx \|\mathbf{x}\|\leq R,  |r|\leq R
\}
$,
$
D_1^{*}= \{ (\mathbf{x},r)\dvtx \|\mathbf{x}\|\leq R,  |r|> R
\}
$, and
$
D_2^{*}= \{ (\mathbf{x},r)\dvtx \|\mathbf{x}\|> R  \}$.
Since $\mathbbmathbf{I}$ and ${\boldsymbol\xi}$ are orthogonal
and $\|\mathbbmathbf{I}\|=1$ we have
$
\| (\mathbf{x},r\mathbbmathbf{I}+s_1(\mathbf
{x}){\boldsymbol\xi} )\|^2
=
\|\mathbf{x}\|^2+r^2+s_1(\mathbf{x})^2\|{\boldsymbol\xi}\|^2\geq\|
\mathbf{x}\|^2+r^2$,
and then
\[
\mathop{\idotsint}_{D_1^{*}}G(\mathbf{x},r) \,\mathrm{d}r\,\mathrm{d}
\mathbf{x} < \mathop{\idotsint}_{\|\mathbf{x}\|\leq R} s_1(
\mathbf{x})^{n_2-1}\,\mathrm{d}\mathbf{x} \times \int_{|r|>R}
\frac{C}{|r|^{n+\delta}}\,\mathrm{d}r <\infty
\]
and
\begin{eqnarray*}
\mathop{\idotsint}_{D_2^{*}} G(\mathbf{x},r)\,\mathrm{d}r\,\mathrm{d}
\mathbf{x}&<& \mathop{\idotsint}_{\|\mathbf{x}\|> R} \int_{-\infty}^{\infty}
\frac{s_1(\mathbf{x})^{n_2-1}}{ (\|\mathbf{x}\|^2+r^2
)^{\trup{(n+\delta)}{2}}}\,\mathrm{d}r \,\mathrm{d}\mathbf{x}
\\
&<&\mathop{\idotsint}_{\|\mathbf{x}\|> R} \frac{s_1(\mathbf{x})^{n_2-1}}{\|\mathbf{x}\|^{n-1+\delta}}\, \mathrm{d}\mathbf
{x} \times \int_{-\infty}^{\infty} \frac{1}{ (1+r^2 )^{n/2}}\,
\mathrm{d}r <\infty,
\end{eqnarray*}
where the multidimensional integral in the last inequality is computed
by means of passing to spherical coordinates.
\end{pf}

Note that in the i.i.d. case the condition
(\ref{Fast_Decay_Criterium_Bound}) is equivalent to the existence of the
$n-1+\delta$ moment of the marginal density $h(x)$. For the Student
one-sample $t$-test, however, the criterium of
Lemma~\ref{Fast_Decay_Criterium} is ``too strict'', see below.

%
\begin{Definition}\label{defA.5}
Multivariate density $g(\mathbf{x})$ has the
asymptotic monotonicity property if there exists a constant $M$ such
that for any $1\leq i \leq n$ and any constants $c_j$, $j\neq i$, the
function $f(x)=g(c_1,\ldots,c_{i-1},x,c_{i+1},\ldots,c_n)$ is
monotone on
$[M,\infty)$.
\end{Definition}

%
\begin{lemma}
\label{Asmptotic_Monotonicity_Criterium}
If $K_g$ is finite and $g(\mathbf{x})$ is bounded and has the
asymptotic monotonicity property, then the assumption
(\ref{OST_Condition_On_Density}) holds.
\end{lemma}
\begin{pf}
Setting $\varepsilon$ equal to $(2\sqrt{n})^{-1}$ and using
asymptotic monotonicity property
we get that the integral in (\ref{OST_Condition_On_Density}) is
bounded by
\[
\int_{0}^{2M\sqrt{n}} r^{n-1} \sup
_{\|{\boldsymbol\xi}\|<\trup{1}{2\sqrt{n}}} g \bigl( r ( \mathbbmathbf{I} + {\boldsymbol\xi} )
\bigr) \,\mathrm{d}r + \int_{2M\sqrt{n}}^{\infty}
r^{n-1} g \biggl( r \frac{
\mathbbmathbf{I}
}{2} \biggr) \,\mathrm{d}r < \infty.
\]
The first summand is finite owing to the boundness of $g$ and the
finiteness of the second summand is
equivalent to the finiteness of $K_g$.
\end{pf}

Asymptotic monotonicity and finiteness of $K_g$ are very mild
constraints. For the i.i.d. case of the Student one-sample $t$-test,
for example,
Lemma~\ref{Asmptotic_Monotonicity_Criterium} implies that the
statement of Theorem~\ref{OST_First_Order_Theorem} holds for any
continuous marginal density $h(\mathbf{x})$ that has monotone tails
and such that $K_g<\infty$, and the latter assumption is weaker than
the assumption of existence of the first moment and holds even for such
heavy tailed densities as Cauchy.

Unfortunately there is no asymptotic monotonicity criterium analogue
for the case of the Student two-sample $t$- and $F$-statistics, and the
constant $K_g$ in (\ref{TST_Def_Kg}) and (\ref{F_Def_Kg}) may be
infinite for some heavy-tailed densities, cf. Bradley \cite{Bradley1952a}.

Finally, in the proofs of Theorems~\ref{OST_First_Order_Theorem},
\ref{TST_First_Order_Theorem} and~\ref{FFirstOrderTheorem} one may have
used the ``almost everywhere'' version of the dominated convergence theorem.
For the Student one-sample $t$-statistic the assumption of continuity
of $g$
can be replaced by the assumption that
$g(\mathbf{x})$ is continuous function of $\mathbf{x}\in\mathbb
{R}^n$ a.e. on the set of points
$
\mathbf{x}=r\mathbbmathbf{I}$, $r>0$,
for the Student two-sample $t$-statistic -- on the set of points
$
\mathbf{x}= r
(
\cos(\omega-\omega_0)\mathbbmathbf{I}_{1}
+
\sin(\omega-\omega_0)\mathbbmathbf{I}_{2}
+
\mathbf{z}
)$,
where $r>0$ and $\omega\in[-\uppi/2,  \uppi/2]$, and for the
$F$-statistic -- on the set of points
$
\mathbf{x}=\mathbb{R}^{n_1}\times r\mathbbmathbf{I}$, $r\in
\mathbb{R}$.
Here a.e. means almost everywhere with respect to the Lebesque measure
induced by the
measure of the linear space $L$ in (\ref{OST_Condition_On_Density}),
(\ref{TST_Condition_On_Density}) and (\ref{F_Condition_On_Density}).\

\section{Figures}
\label{SectionFiguresandTables}

\setcounter{figure}{1}
\renewcommand{\thefigure}{\arabic{figure}}
%
\begin{figure}[h!]

\includegraphics{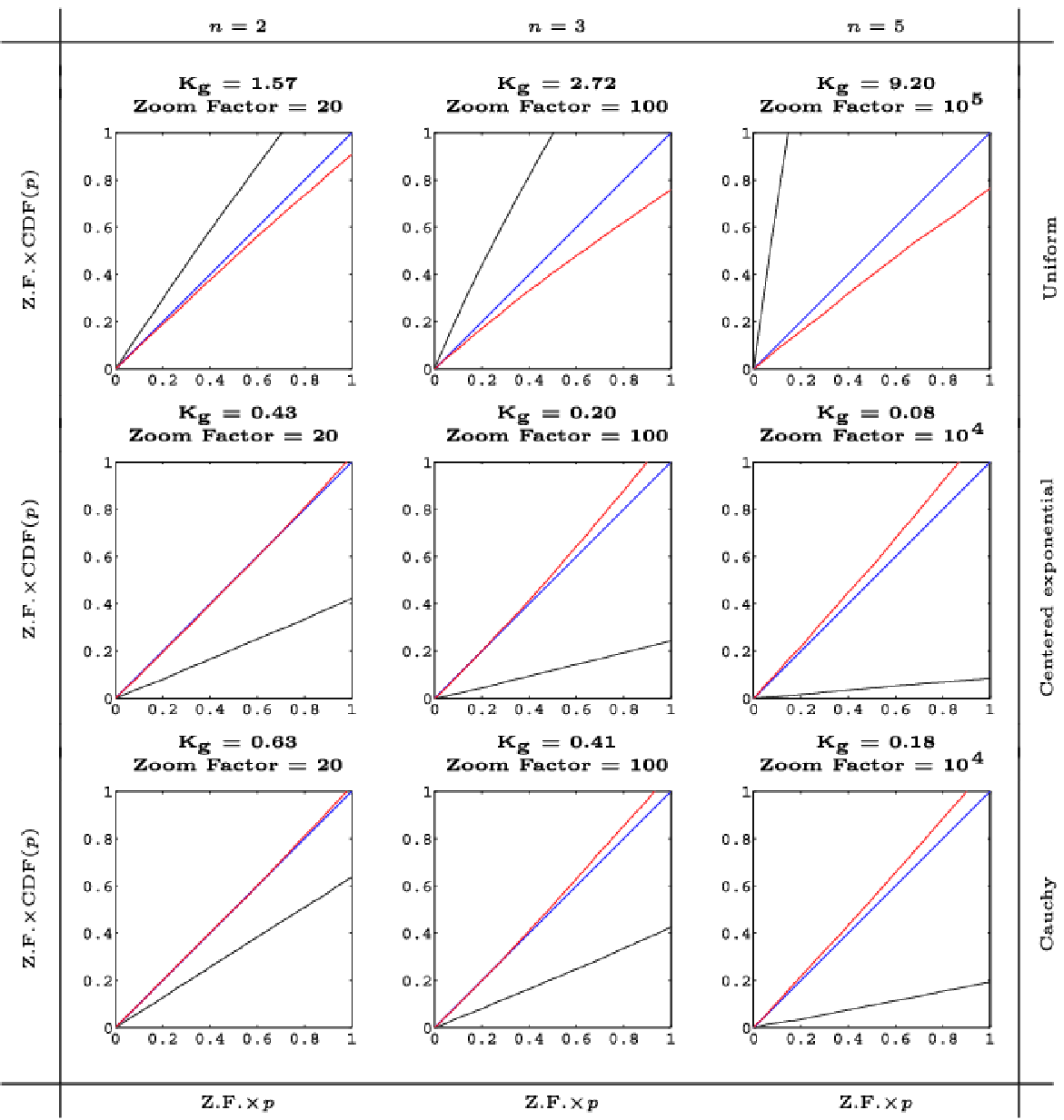}

%
\caption{The eCDF of the $p$-values for the Student one-sample $t$-test.
The empirical CDFs of the raw and corrected $p$-values $p^R$ and $p^C$
are shown in black and red accordingly.
The top, middle and bottom rows correspond to the $\operatorname
{Uniform}(-1,1)$,
\textit{Centered exponential} and \textit{Cauchy} densities,
and left, middle and right columns correspond to sample sizes $n=2$,
$n=3$ and $n=5$.
The blue diagonal line is the theoretical uniform distribution. The
axes are scaled according to the Zoom Factor (Z.F.) parameter $r$
in the title of the graphs.}\label{fig:SimulationOST}\vspace*{51pt}
\end{figure}

%
%
\begin{figure}[t!]

\includegraphics{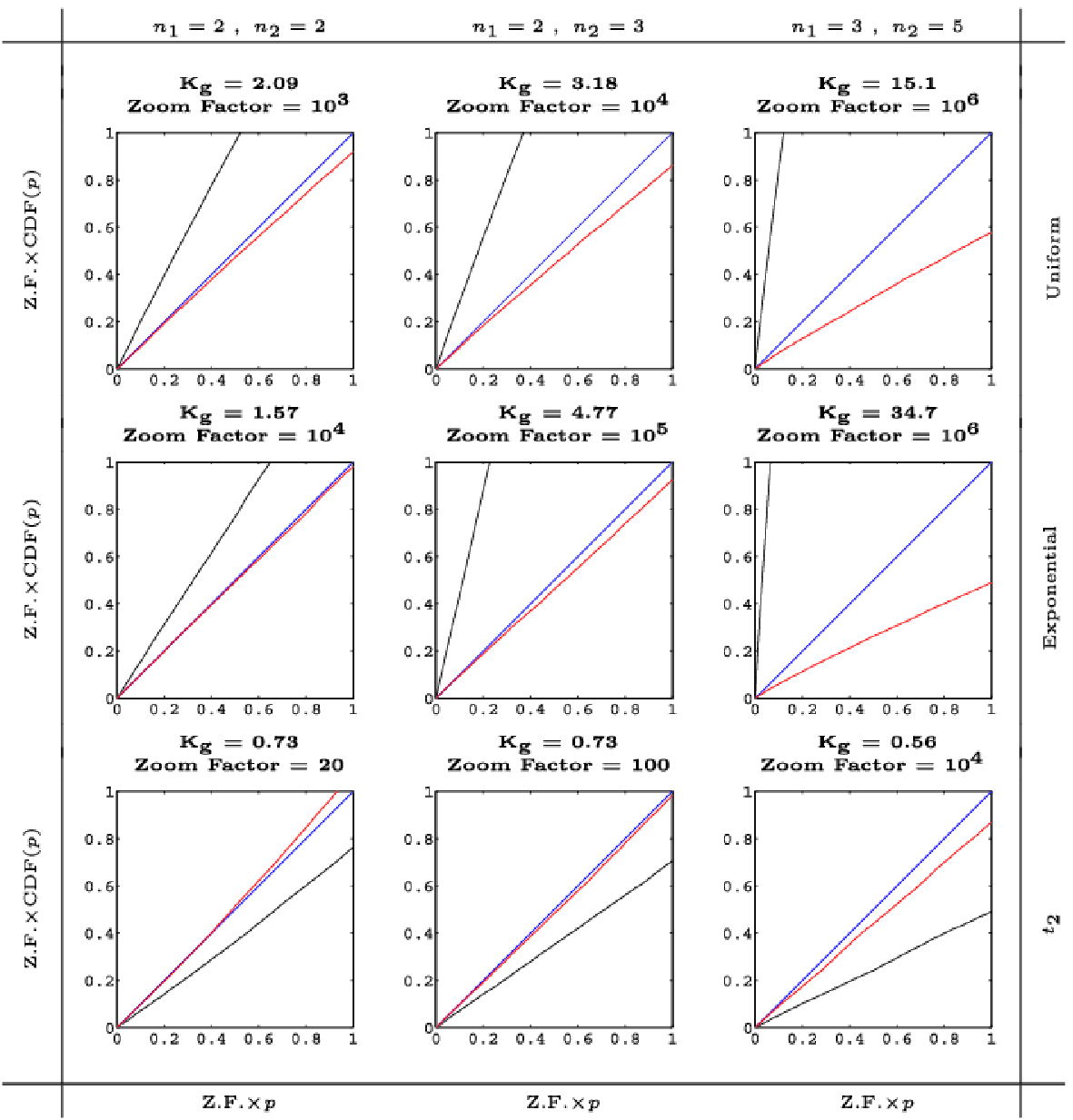}

%
%
%
\caption{The eCDF of the $p$-values for the Student two-sample $t$-test.
The empirical CDFs of the raw and corrected $p$-values $p^R$ and $p^C$
are shown in black and red accordingly.
The top, middle and bottom rows correspond to the $\operatorname
{Uniform}(-1,1)$,
\textit{Exponential} and $t_2$ densities,
and left, middle and right columns correspond to sample sizes
($n_1=2$, $n_2=2$), ($n_1=2$, $n_2=3$), and ($n_1=3$, $n_2=5$).
The blue diagonal line is the theoretical uniform distribution.
The axes are scaled according to the Zoom Factor (Z.F.) parameter $r$
in the title of the graphs.}\label{fig:SimulationTST}\vspace*{60pt}
\end{figure}

%
%
\begin{figure}[t!]

\includegraphics{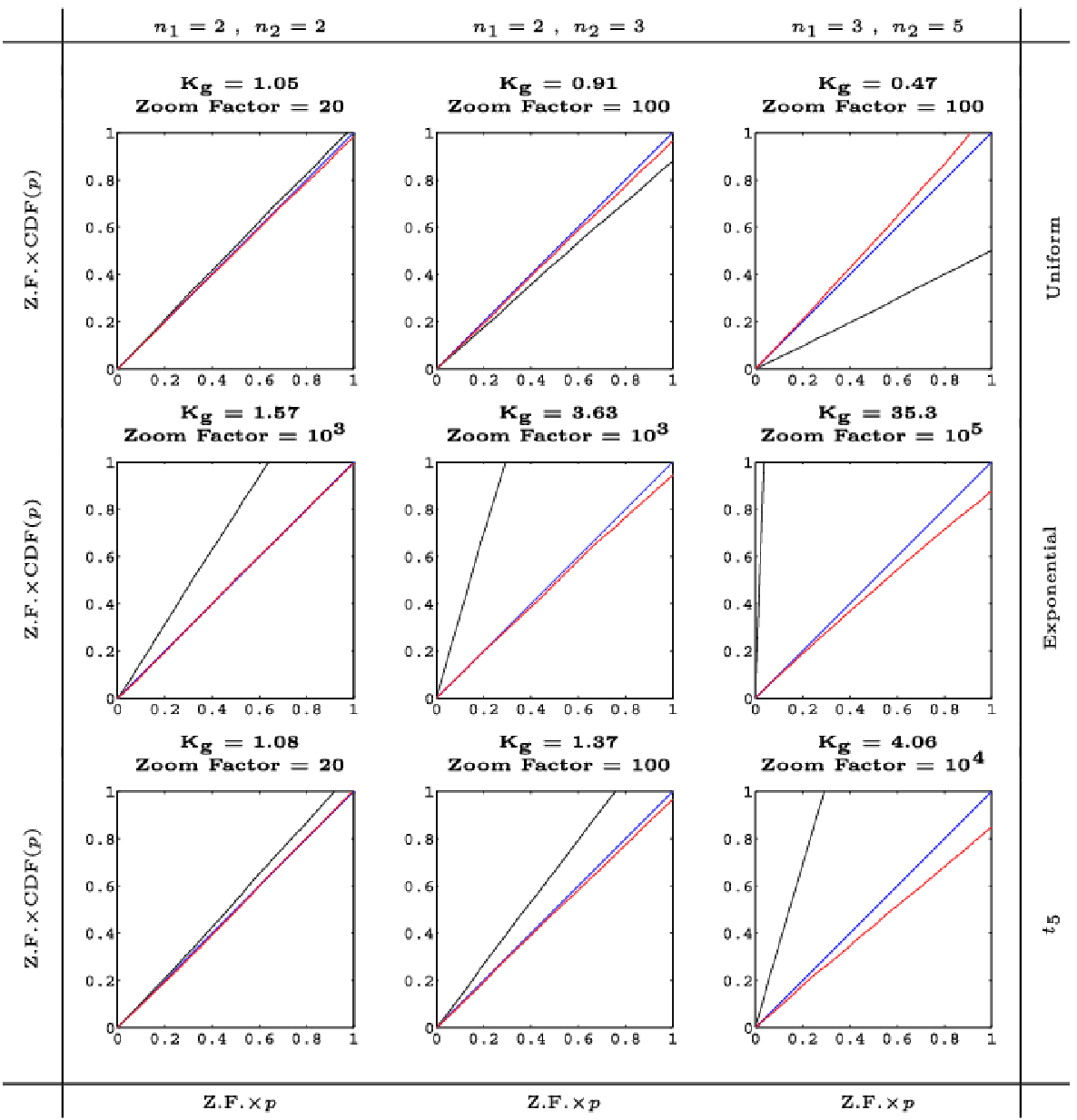}

%
\caption{The eCDF of the $p$-values for the $F$-test (equality of variances).
The empirical CDFs of the raw and corrected $p$-values $p^R$ and $p^C$
are shown in black and red accordingly.
The top, middle and bottom rows correspond to the $\operatorname
{Uniform}(-1,1)$,
\textit{Exponential} and $t_5$ densities,
and left, middle and right columns correspond to sample sizes
($n_1=2$, $n_2=2$), ($n_1=2$, $n_2=3$), and ($n_1=3$, $n_2=5$).
The blue diagonal line is the theoretical uniform distribution.
The axes are scaled according to the Zoom Factor (Z.F.) parameter $r$
in the title of the graphs.}\label{fig:SimulationF}\vspace*{60pt}
\end{figure}

%
%
\begin{figure}[t!]

\includegraphics{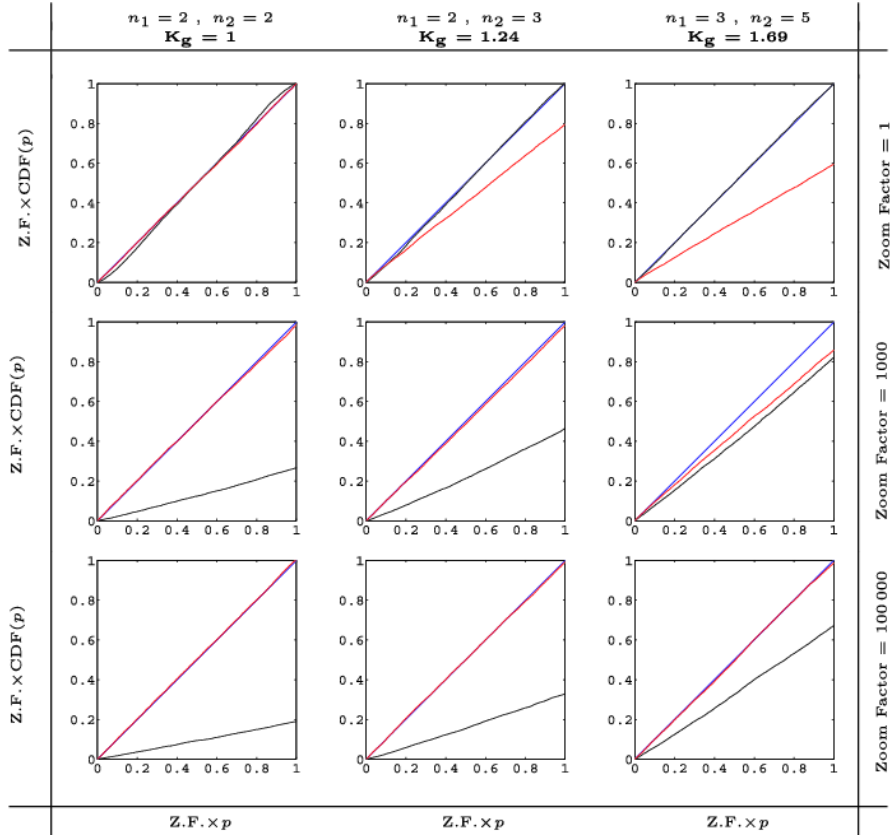}

%
%
%
%
\caption{The distribution tails of the $p$-values for the Welch test.
The empirical CDFs of the raw (Welch--Satterthwaite) and corrected
$p$-values $p^R$ and $p^{\mathit{WS}}$ for the \textit{Standard Normal} density
are shown in black and red accordingly.
The top, middle and bottom rows correspond to the different values of
the Zoom Factor (Z.F.) parameter $r$ shown on the right, and
the axes are scaled accordingly. The left, middle and right columns
correspond to sample sizes ($n_1=2$, $n_2=2$), ($n_1=2$, $n_2=3$), and
($n_1=3$, $n_2=5$). The blue diagonal line is the theoretical uniform
distribution.}\label{fig:SimulationWELCH}\vspace*{90pt}
\end{figure}

%
%
\begin{figure}[t!]

\includegraphics{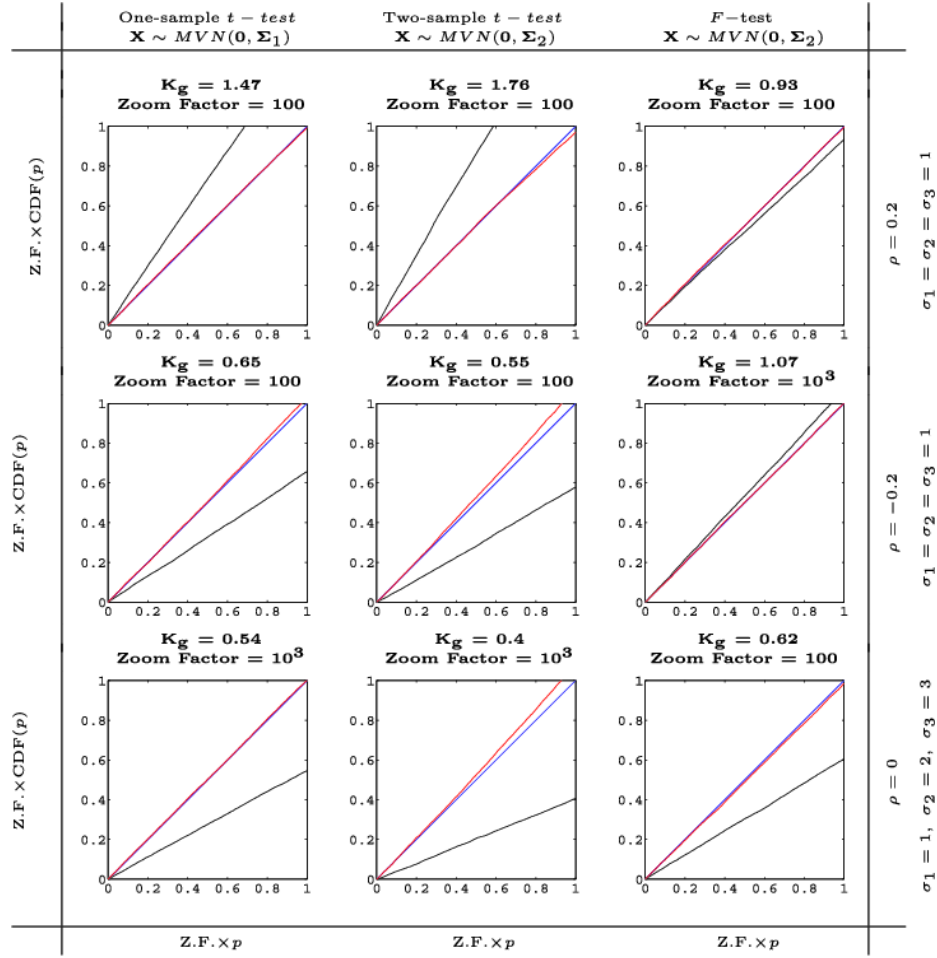}

%
%
%
\caption{The effect of dependency and non-homogeneity of data on
$\mathbf{P}({T>u})$ as $u\to\infty$.
The empirical CDF of raw (black) and corrected (red) $p$-values.
Analogue of Figures~\protect\ref{fig:SimulationOST},
\protect\ref{fig:SimulationTST}
and \protect\ref{fig:SimulationF} for dependent (top row -- positively
correlated observations; middle row -- negatively correlated
observations) and non-homogeneous (bottom row, unequal variances) data.
Multivariate normal case with covariance\vspace*{12pt}
matrices\newline
\mbox{\qquad}$\displaystyle\boldsymbol{\Sigma}_1=
\lleft(
\protect\matrix{
\sigma_1^2 & \rho\sigma_1\sigma_2 & 0
\cr
\rho\sigma_1\sigma_2 & \sigma_2^2 & \rho\sigma_2\sigma_3
\cr
0 & \rho\sigma_2\sigma_3 & \sigma_3^2
}
\rright),
\qquad
\displaystyle\boldsymbol\Sigma_2=
\lleft(
\protect\matrix{
\sigma_1^2 & \rho\sigma_1\sigma_2 & 0 & 0 & 0
\cr
\rho\sigma_1\sigma_2 & \sigma_2^2 & 0 & 0 & 0
\cr
0 & 0 & \sigma_1^2 & \rho\sigma_1\sigma_2 & 0
\cr
0 & 0 & \rho\sigma_1\sigma_2 & \sigma_2^2 & \rho\sigma_2\sigma_3
\cr
0 & 0 & 0 & \rho\sigma_2\sigma_3 & \sigma_3^2
}
\rright)$.}\label{fig:SimulationMVN}
\end{figure}

\end{appendix}

\section*{Acknowledgements}

The author thanks Holger Rootz\'{e}n for assistance and
fruitful discussions.
He also thanks the associate editor for extremely helpful comments which
led to improvement of the presentation of the present paper.

\begin{supplement}\label{SectionSupplementaryMaterials}
\stitle{MATLAB, Wolfram Mathematica scripts, other materials\\}
\slink[doi]{10.3150/13-BEJ552SUPP} 
\sdatatype{.zip}
\sfilename{bej552\_supp.zip}
\sdescription{\mbox{}\\[12pt]
\textbf{MATLAB scripts}.
[\textit{OST}/\textit{TST}/\textit{WELCH}/\textit{F}]$+$\textit
{ComputeKg.m} -- compute $K_g$ for the Student one- and two-sample
$t$-, Welch, and $F$-statistics using adaptive Simpson or Lobatto
quadratures. Here $g$ is an arbitrary multivariate
density.\footnote{For the $F$-statistic we use Monte Carlo integration.}
[\textit{TST}/\textit{WELCH}/\textit{F}]$+\break$\textit{ComputeKgIS}$+$.m -- the same as above but for the case where samples are
independent.$^{2}$
$[$\textit{OST}/\textit{TST}/\textit{WELCH}/\textit{F}$]+$\textit
{ComputeKgIID}$+$.m -- the same as above but assuming that the samples
consist of i.i.d. random variables.\footnote{For the $F$-statistic and
$n_1>3$ we use Monte Carlo integration.}
\textit{RunSimulation}$+$[\textit{IID/MVN}]$+$.m -- perform simulation
study for i.i.d. and dependent/non-homogeneous cases, see Section~\ref
{SectionSimulation} and \hyperref[SectionFiguresandTables]{Appendix B}.\\[12pt]
\textbf{Wolfram Mathematica scripts}.
[\textit{OST}/\textit{TST}/\textit{WELCH}/\textit{F}]$+$\textit
{ComputeKg.nb} -- compute the exact expression for $K_g$ for an
arbitrary multivariate density
$g$ and given sample size(s). We include a number of examples, such as
evaluation of $K_g$ for the zero-mean Gaussian case with an arbitrary
covariance matrix $\boldsymbol\Sigma$; the ``unequal variances'' case
for the Student two-sample $t$- and Welch statistics; and evaluation of
$K_g$ for the densities considered in the simulation study.
\mbox{\textit{OSTComputeKgIID.nb}} -- verifies the constants in Table~\ref
{tab:KgOSTIID} for the i.i.d. case of the Student one-sample $t$-statistic.
\textit{TSTExactPDF.nb} and \textit{WELCHExactPDF.nb} -- the exact
distribution for the Student two-sample $t$- and Welch statistics for
odd sample sizes, see Ray and Pitman \cite{Ray1961}.\\[12pt]
\textbf{Other materials}.
\textit{Supplementary-Materials.pdf} -- Remarks on Theorem~\ref{IntroductionMainTheorem} and its application to real data; extended
version of the literature review; comparison of the result of Theorem~\ref{IntroductionMainTheorem} with the exact distribution of the
Welch statistic; proof of Theorem~\ref{FFirstOrderTheorem}.}
\end{supplement}



%
%

\printhistory

\end{document}